\documentclass[pdflatex,sn-basic,Numbered]{sn-jnl}
 
\usepackage{bm}

\usepackage{graphicx}%
\usepackage{multirow}%
\usepackage{amsmath,amssymb,amsfonts}%
\usepackage{amsthm}%
\usepackage{mathrsfs}%
\usepackage[title]{appendix}%
\usepackage{xcolor}%
\usepackage{textcomp}%
\usepackage{manyfoot}%
\usepackage{booktabs}%
\usepackage{algorithm}%
\usepackage{algpseudocode}%
\usepackage{listings}%
\usepackage{cleveref}
\usepackage{verbatim}
\usepackage{caption}
\usepackage{subcaption}
\usepackage{float}

\def\N{\mathbb{N}}

\newcommand{\bsgamma}{\boldsymbol{\gamma}}
\newcommand{\bszero}{\boldsymbol{0}}
\newcommand{\bsDelta}{\boldsymbol{\Delta}}

\newcommand{\bsk}{\boldsymbol{k}}
\newcommand{\bsell}{\boldsymbol{\ell}}

\newcommand{\bsa}{\boldsymbol{a}}

\newcommand{\bsg}{\boldsymbol{g}}
\newcommand{\bsx}{{\boldsymbol{x}}}
\newcommand{\bsy}{{\boldsymbol{y}}}

\newcommand{\EE}{\mathbb{E}}
\newcommand{\R}{\mathbb{R}}
\newcommand{\ZZ}{\mathbb{Z}}

\numberwithin{equation}{section}

\newcommand{\rd}{\, \mathrm{d}}
\newcommand{\ran}{\mathrm{ran}}
\newcommand{\wor}{\mathrm{wor}}

\allowdisplaybreaks

\theoremstyle{thmstyleone}%
\newtheorem{theorem}{Theorem}
\theoremstyle{thmstyletwo}%
\newtheorem{remark}{Remark}%

\theoremstyle{thmstylethree}%
\newtheorem{lemma}{Lemma}%

\begin{document}

\title[A lattice algorithm with multiple shifts for function approximation in Korobov spaces]{A lattice algorithm with multiple shifts for function approximation in Korobov spaces}


\author[1]{\fnm{Mou} \sur{Cai}}\email{caimoumou@g.ecc.u-tokyo.ac.jp}
\equalcont{These authors contributed equally to this work.}

\author*[2]{\fnm{Josef} \sur{Dick}}\email{josef.dick@unsw.edu.au}
\equalcont{These authors contributed equally to this work.}

\author[1]{\fnm{Takashi} \sur{Goda}}\email{goda@frcer.t.u-tokyo.ac.jp}
\equalcont{These authors contributed equally to this work.}

\affil*[1]{\orgdiv{Graduate School of Engineering}, \orgname{The University of Tokyo}, \orgaddress{\street{7-3-1 Hongo, Bunkyo-ku}, \city{Tokyo}, \postcode{113-8656}, 
\country{Japan}}}

\affil[2]{\orgdiv{School of Mathematics and Statistics}, \orgname{UNSW Sydney}, \orgaddress{
\city{Sydney}, \postcode{2052}, \state{NSW}, \country{Australia}}}



\abstract{In this paper, we propose a novel algorithm for function approximation in a weighted Korobov space based on shifted rank-1 lattice rules. To mitigate aliasing errors inherent in lattice-based Fourier coefficient estimation, we employ $\mathcal{O}((\log N)^{2d-1})$ shifted copies of a single rank-1 lattice and recover each Fourier coefficient via a least-squares procedure. Writing $p$ for the total number of function evaluations, we show that the resulting approximation achieves the optimal convergence rate for the $L_{\infty}$-approximation error in the worst-case setting, namely $\mathcal{O}(p^{-\alpha+1/2+\varepsilon})$ for arbitrarily small $\varepsilon>0$. Moreover, by incorporating random shifts, the algorithm attains the optimal rate for the $L_{2}$-approximation error in the randomized setting, which is $\mathcal{O}(p^{-\alpha+\varepsilon})$. Numerical experiments illustrate the practical performance of the algorithms and the qualitative behavior predicted by the theoretical analysis.}




\maketitle

\section{Introduction}

We study approximation of smooth periodic functions defined over the $d$-dimensional unit cube $[0,1)^d$ based on pointwise function evaluations. For an $N$-element point set $P=\{\bsx_1,\ldots,\bsx_N\}\subset [0,1)^d$, a general approximation scheme takes the form
\begin{align}\label{eq:form_approximation}
    f(\bsx)\approx \mathcal{A}(f)(\bsx)=\sum_{i=1}^{N}c_i(\bsx)f(\bsx_i),
\end{align}
where the coefficient functions $c_1,\ldots,c_N: [0,1)^d\to \mathbb{C}$ determine the approximation. The error of such an approximation is commonly measured in an $L_p$ norm, namely
\[ \|f-\mathcal{A}(f)\|_{q} = \left( \int_{[0,1)^d}\left|f(\bsx)-\mathcal{A}(f)(\bsx)\right|^{q} \rd \bsx\right)^{1/q},\]
for $q\in [1,\infty)$, with the usual modification for $q=\infty$. 

For a Banach space $F$ with norm $\|\cdot\|_F$, the worst-case error of a deterministic algorithm $\mathcal{A}$ is defined as
\[ \mathrm{e}_{q}^{\wor}(\mathcal{A},F) := \sup_{\substack{f\in F\\ \|f\|_F\le 1}}\|f-\mathcal{A}(f)\|_{q}. \]
When we speak of a randomized algorithm, we refer to a pair consisting of a probability space $(\Omega, \Sigma, \mu)$ and a family of mappings $\mathcal{A} = (\mathcal{A}^\omega)_{\omega \in \Omega}$, where each $\mathcal{A}^\omega$ is a deterministic algorithm of the form \eqref{eq:form_approximation} for fixed $\omega \in \Omega$. Then, for a Banach space $F$, the random-case error of a randomized algorithm $\mathcal{A}$ is defined as
\[ \mathrm{e}_{q}^{\ran}(\mathcal{A},F) := \sup_{\substack{f\in F\\ \|f\|_F\le 1}}\left(\EE_{\omega}\left[ \|f-\mathcal{A}^\omega(f)\|^{q}_{q}\right]\right)^{1/q},\]
where the expectation is taken with respect to the measure $\mu$.

In what follows, we move from this abstract setting to the specific function space of interest, namely the weighted Korobov space $\mathcal{K}_{d,\alpha,\bsgamma}$ with smoothness $\alpha>1/2$ and positive weights $\bsgamma=(\gamma_1,\gamma_2,\ldots)$. As explained in Section~\ref{subsec:space}, the space contains smooth periodic functions over $[0,1)^d$, and every function in the space has an absolutely convergent Fourier series
\[ f(\bsx) = \sum_{\bsk \in \mathbb{Z}^d} \widehat{f}(\bsk)\, \mathrm{e}^{2\pi \mathrm{i} \bsk \cdot \bsx},
\]
where $\widehat{f}(\bsk)$ denotes the $\bsk$-th Fourier coefficient of $f$:
\[ \widehat{f}(\bsk):=\int_{[0,1)^d}f(\bsx)\, \mathrm{e}^{-2\pi \mathrm{i} \bsk \cdot \bsx}\rd \bsx. \]
For an index set $A\subset \mathbb{Z}^d$ and an $N$-element point set $P$, we consider approximating $f\in \mathcal{K}_{d,\alpha,\bsgamma}$ by
\begin{align}\label{eq:Fourier_approximation}
    f(\bsx) & \approx \sum_{\bsk \in A} \widehat{f}(\bsk)\, \mathrm{e}^{2\pi \mathrm{i} \bsk \cdot \bsx} \approx \sum_{\bsk \in A} \left( \sum_{i=1}^{N}w_{i,\bsk}f(\bsx_i)\, \mathrm{e}^{-2\pi \mathrm{i} \bsk \cdot \bsx_i}\right) \mathrm{e}^{2\pi \mathrm{i} \bsk \cdot \bsx}:=\mathcal{A}(f)(\bsx),
\end{align}
with quadrature weights $\{w_{1,\bsk},\ldots,w_{N,\bsk}\}_{\bsk\in A}$.
That is, we first restrict the range of $\bsk\in \mathbb{Z}^d$ to $A$, and then for each $\bsk \in A$, we approximate the Fourier coefficient $\widehat{f}(\bsk)$ by quadrature using function values on $P$.

The key issue inherent in using a single rank-1 lattice point set is \emph{aliasing}. For a single rank-1 lattice rule, one can simplify the computation of the Fourier coefficients in \eqref{eq:Fourier_approximation} as
\begin{align}\label{eq:Aliasing_state}
\widehat{f}(\bsk)\approx \frac{1}{N} \sum_{i=1}^{N} f(\bsx_i)\, \mathrm{e}^{-2\pi \mathrm{i} \bsk \cdot \bsx_i}=\sum_{\bsell \in \mathcal{L}^{\perp}}\widehat{f}(\bsk+\bsell),
\end{align}
where $\mathcal{L}^{\perp}$ denotes the dual lattice of $P$, whose definition will be given in Section~\ref{subsec:lattice}.
If there exist distinct frequencies $\bsk,\bsk'\in A$ and some $\bsell\in \mathcal{L}^{\perp}$ such that $\bsk+\bsell=\bsk'$, then the approximation of the corresponding Fourier coefficients coincides and their approximations cannot be distinguished. To overcome the aliasing problem, the number of sampling points must be at least as large as the number of Fourier coefficients to be approximated. We refer to this requirement as oversampling. In the following, we denote the number of function evaluations divided by the number of Fourier coefficients which we approximate as the oversampling factor.

To overcome this aliasing issue, several approaches have been proposed. In \cite{KV19}, K\"{a}mmerer and Volkmer introduced the use of multiple rank-1 lattice point sets as $P$, and showed that, for $q=\infty$, the resulting algorithm of the form \eqref{eq:Fourier_approximation} achieves the optimal rate $\mathcal{O}(N^{-\alpha+1/2+\varepsilon})$ of the worst-case error with an oversampling factor of $\mathcal{O}(\log |A|)$. In \cite{PKG25}, Pan et al.\ considered choosing multiple shifted rank-1 lattice point sets independently and randomly, and then estimating each Fourier coefficient $\widehat{f}(\bsk)$ for $\bsk\in A$ by taking the median over the independent estimates. They obtained, for $q=2$, an error bound of $\mathcal{O}(N^{-\alpha+\varepsilon})$ for any fixed function $f\in \mathcal{K}_{d,\alpha,\bsgamma}$, with high probability, and the optimal rate $\mathcal{O}(N^{-\alpha+\varepsilon})$ for the random-case error, with an oversampling factor of roughly $\mathcal{O}((\log |A|)^2)$.
Another line of work to mitigate aliasing was recently introduced in \cite{BGKS25}, where a single rank-1 lattice point set of substantially larger size is first generated, and then $N$ subsamples are randomly selected. They showed that, for $q=2$, the optimal rate $\mathcal{O}(N^{-\alpha+\varepsilon})$ for the worst-case error is achieved with high probability. The oversampling factor in that paper is $\mathcal{O}(\log |A|)$.

The aim of this paper is to introduce a new approach to mitigating aliasing errors inherent in lattice-based Fourier coefficient estimation. 
The key idea of our approach is to partition the index set $A$ into disjoint subsets in such a way that, within each subset, all frequencies are aliased under a single rank-1 lattice rule.  For all subsets together, we use one common family of at most $\mathcal{O}((\log N)^{2d-1})$ shifted copies of the rank-1 lattice, from which the Fourier coefficients are recovered through separate least-squares problems on the individual subsets.
The main difference from the multiple rank-1 lattice algorithm in \cite{KV19} is that, instead of taking the union of several rank-1 lattice point sets, our method relies on a single rank-1 lattice point set. The aliasing issue is then addressed by applying altogether at most $\mathcal{O}((\log N)^{2d-1})$ shifted copies, so that the total number of function evaluations is $p=\mathcal{O}(N(\log N)^{2d-1})$, and subsequently performing a least-squares procedure to estimate the Fourier coefficients.
We prove that the proposed algorithm achieves the optimal rate $\mathcal{O}(p^{-\alpha+1/2+\varepsilon})$ for the worst-case $L_\infty$-approximation error (Theorem~\ref{thm_det}). Furthermore, by incorporating random shifts, we show that the resulting randomized algorithm attains the optimal rate $\mathcal{O}(p^{-\alpha+\varepsilon})$ for the randomized $L_2$-approximation error (Theorem~\ref{thm:L_2_rce}). Although the theoretical sampling overhead relative to the underlying lattice size is substantially larger than that of the methods discussed above, our numerical experiments indicate that such a large overhead is not required in practice; see Section~\ref{sec:numerical_experiments}.

The remainder of the paper is organized as follows. Section~\ref{sec:notation_background} reviews the necessary background, and Section~\ref{sec:algorithm} presents the algorithms and proves that suitable generating vectors and shifts can be found with high probability. Section~\ref{sec:error_analysis} establishes optimal convergence rates for both the deterministic $L_\infty$ error and the randomized $L_2$ error. We conclude this paper with numerical experiments in Section~\ref{sec:numerical_experiments}, which illustrate the practical performance of the proposed algorithms and the qualitative behavior predicted by the theoretical analysis.

\subsection{Literature review}

Originally, in \cite{KSW06}, Kuo, Sloan and Wo\'zniakowski studied approximations of the form \eqref{eq:Fourier_approximation}, with $P$ being a single rank-1 lattice point set and $w_{i,\bsk}=1/N$ for all $i$ and $\bsk$. They proved that, for $q=2$, the resulting algorithm achieves a worst-case error of $\mathcal{O}(N^{-\alpha/2+\varepsilon})$ for arbitrarily small $\varepsilon>0$. Later, Byrenheid, K\"ammerer, Ullrich, Volkmer showed in \cite{BKUV17} that,  as long as a single rank-1 lattice point set is used as $P$, any general algorithm of the form \eqref{eq:form_approximation} cannot attain a convergence rate better than $\mathcal{O}(N^{-\alpha/2})$, which is only half of the best possible rate $\mathcal{O}(N^{-\alpha +\varepsilon})$. By employing randomized single rank-1 lattice point sets, Cai, Goda and Kazashi~\cite{CGK24} recently proved that, for $q=2$, the randomized lattice-based algorithm achieves a random-case error of $\mathcal{O}(N^{-\alpha(2\alpha+1)/(4\alpha+1)+\varepsilon})$, which improves upon $\mathcal{O}(N^{-\alpha/2})$ for any $\alpha>1/2$, yet remains far from the optimal rate $\mathcal{O}(N^{-\alpha + \varepsilon})$. The case $q=\infty$ has also been studied, for instance in \cite{KWW09,ZKH09}, where a single rank-1 lattice point set is used as $P$, and a worst-case error of $\mathcal{O}(N^{-\alpha(2\alpha-1)/(4\alpha-1)+\varepsilon})$ was obtained. Again, this convergence rate is far from the optimal rate $\mathcal{O}(N^{-\alpha+1/2+ \varepsilon})$ achievable in the case $q=\infty$. Fast CBC constructions underpin most practical implementations \cite{NC06}. See also the modern textbook treatment which consolidates these ideas and their extensions \cite{DKP22}.

To circumvent aliasing on a single lattice, one may combine a few \emph{different} lattices. K{\"a}mmerer and Volkmer proved that sampling along multiple rank-1 lattices removes aliasing and achieves the \emph{optimal} worst-case rate $\mathcal O(N^{-\alpha+1/2+\varepsilon})$ for $L_\infty$-approximation in Korobov spaces \cite{KV19}. We also refer to a recent work \cite{CG26} regarding the error and tractability analyses in the weighted Korobov spaces. This line builds on a reconstruction theory guaranteeing that hyperbolic-cross trigonometric polynomials can be stably recovered from rank-1 lattice samples \cite{Kammerer2013SJNA} and on single-lattice approximation bounds for hyperbolic-cross index sets \cite{KammererPottsVolkmer2015JC}. The algorithmic foundations of multiple rank-1 lattices for exact recovery (and their use in sparse FFT-style methods) are developed in \cite{Kammerer2018JFAA}, and a deterministic near-optimal construction was given in \cite{GrossIwenKammererVolkmer2021}.

Randomization offers another route to de-aliasing while keeping a single generating vector. For $L_2$-approximation, Cai, Goda and Kazashi analyzed a randomized single-lattice method and proved rates strictly better than the deterministic single-lattice barrier, though still below optimal \cite{CGK24}. Pan, Kritzer and Goda subsequently introduced a robust \emph{median-of-means} estimator that evaluates the same lattice under many independent random shifts; they established high-probability and randomized $L_2$ error bounds at the optimal rate $\mathcal O(N^{-\alpha+\varepsilon})$, with a follow-up “universal’’ result that eliminates tuning to the smoothness \cite{PKG25,PanGodaKritzer2025Universal}. Random shifting itself is classical in randomized QMC, beginning with Cranley–Patterson and refined in later distributional analyses; these provide the probabilistic underpinning for shift-based variance/error control \cite{CranleyPatterson1976,LEcuyerLemieux2010}.

A complementary strategy is to start from one \emph{very large} lattice and then select $N$ points uniformly at random. Bartel, Gilbert, Kuo and Sloan proved that such minimal random subsamples achieve the optimal $\mathcal O(N^{-\alpha+\varepsilon})$ rate for worst-case $L_2$-approximation with high probability, giving another single-lattice path to optimality \cite{BGKS25}.

Least-squares reconstruction from lattice samples sits naturally inside a general reproducing kernel Hilbert space (RKHS) framework and also extends to nonperiodic settings via tent and cosine transforms. Kuo, Migliorati, Nobile and Nuyens provide a unified view of integration, reconstruction and approximation from lattice data, together with efficient CBC strategies tailored to these tasks \cite{KuoMiglioratiNobileNuyens2021}. Finally, the optimal exponents quoted above match the best-possible rates known from hyperbolic-cross approximation and sampling-number theory for spaces with dominating mixed smoothness; see, e.g., the survey \cite{Dung2016HC} and recent sampling-number advances \cite{JahnUllrichVoigtlaender2022}.

Within this landscape, our paper addresses aliasing \emph{within the fibers} of a single rank-1 lattice by using a common family of at most $\mathcal O((\log N)^{2d-1})$ carefully chosen shifted copies and recovering the Fourier coefficients through small least-squares solves for the individual fibers. This preserves the sampling simplicity of a single lattice (as in randomized approaches) while matching the optimal rates previously known from multiple-lattice schemes in $L_\infty$ and from randomized analyses in $L_2$.

\section{Notation and background}\label{sec:notation_background}

Throughout the paper we denote vectors in $\ZZ^d$ or $\R^d$ by boldface letters such as $\bsk=(k_1,\ldots,k_d)$. For $\bsx,\bsy\in \R^d$, we denote their standard inner product by $\bsx\cdot \bsy = x_1 y_1 + \cdots + x_d y_d$. For a subset $u\subseteq \{1,\ldots,d\}$, $\bsx_u$ denotes the vector $(x_j)_{j\in u}$. For a complex matrix $B\in \mathbb{C}^{m\times n}$, we denote its Hermitian transpose by $B^H = \overline{B}^\top$. 

In the following, we first define the weighted Korobov space $\mathcal{K}_{d,\alpha,\bsgamma}$ and introduce rank-1 lattice point sets and some of their properties.

\subsection{Weighted Korobov space}\label{subsec:space}

For a smoothness parameter $\alpha > 1/2$ and positive weights $\bsgamma=(\gamma_1,\gamma_2,\ldots)$, we define
\[
 r_{\alpha,\bsgamma}(\bsk) := \prod_{j\in u(\bsk)} \frac{|k_j|^\alpha}{\gamma_j}, 
\]
for $\bsk\in \ZZ^d$, where we set $u(\bsk) = \{j \in \{1, 2, \ldots, d\}: k_j \neq 0\}$, and the empty product is set to $1$. 

Let $f:[0,1)^d \to \mathbb{C}$ be a $1$-periodic function given by an absolutely convergent Fourier series
\[
 f(\bsx) = \sum_{\bsk \in \ZZ^d} \widehat{f}(\bsk)\, \mathrm{e}^{2\pi \mathrm{i} \bsk \cdot \bsx},
\]
where $\widehat{f}(\bsk)$ denotes the $\bsk$-th Fourier coefficient.
We define the weighted Korobov space $\mathcal{K}_{d,\alpha,\bsgamma}$ with smoothness parameter $\alpha$ as the set of $1$-periodic functions whose Fourier coefficients satisfy
\begin{align}\label{eq:kor_norm}
\|f\|_{d,\alpha, \bsgamma} := \left( \sum_{\bsk \in \mathbb{Z}^d} |\widehat{f}(\bsk)|^2\, r_{\alpha,\bsgamma}^2(\bsk) \right)^{1/2} < \infty.
\end{align}

Note that $\mathcal{K}_{d,\alpha,\bsgamma}$ is a reproducing kernel Hilbert space, with reproducing kernel
\[
K_{d,\alpha,\bsgamma}(\bsx,\bsy):=\sum_{\bsk \in \ZZ^d}\frac{\mathrm{e}^{2\pi \mathrm{i} \bsk \cdot (\bsx-\bsy)}}{r^2_{\alpha,\bsgamma}(\bsk)}
\]
and inner product
\[
\langle f,g\rangle_{d,\alpha,\bsgamma}=\sum_{\bsk \in \mathbb{Z}^d} \widehat{f}(\bsk)\overline{\widehat{g}(\bsk)}\, r^2_{\alpha,\bsgamma}(\bsk).
\]
Furthermore, if $\alpha>1/2$ is an integer, the norm $\|f\|_{d,\alpha, \bsgamma}$ can alternatively be expressed in terms of  (square-integrable) mixed partial derivatives of $f$ up to order $\alpha$ in each variable; see for instance \cite[Proposition~2.17]{DKP22}.

Regarding the index set $A$ in \eqref{eq:Fourier_approximation}, we consider the following special form: for a positive real number $M > 0$, let
\begin{align}\label{eq:Def_hyperbolic_set}
A:=A_{\alpha, \bsgamma, M} = \{\bsk \in \mathbb{Z}^d: r_{\alpha,\bsgamma}(\bsk) < M \},
\end{align}
where $M$ will be chosen appropriately later, depending on $N,d,\alpha,\bsgamma$. 
The following lemma, adapted from \cite[Lemma~1]{KSW06}, provides an upper bound on the number of elements in $A_{\alpha, \bsgamma, M}$.

\begin{lemma}\label{lem:hyp_cross_size}
For $d\in \N$, $\alpha>1/2$, positive weights $\bsgamma=(\gamma_1,\gamma_2,\ldots)$, and $M>0$, let $A_{\alpha,\bsgamma,M}$ be defined by \eqref{eq:Def_hyperbolic_set}. Then, for any $\lambda>1/ \alpha$, we have 
\[ 
|A_{\alpha,\bsgamma, M}|\leq M^{ \lambda} \prod_{j=1}^d \left(1+2\gamma_j^{\lambda}\zeta(\alpha\lambda)\right),
\]
where $\zeta(s) = \sum_{k=1}^\infty k^{-s}$ denotes the Riemann zeta function.
\end{lemma}

If we choose $M$ as the infimum of all $M^\ast > 0$ such that $|A_{\alpha, \bsgamma, M^\ast}| > N$, that is
\begin{equation}\label{eq:def_M}
M := M(\alpha,\bsgamma,N)= \inf\{ M^\ast \in \mathbb{R}_{>0}: |A_{\alpha, \bsgamma, M^\ast}| > N\},
\end{equation}
then $|A_{\alpha, \bsgamma, M}| \le N$. 
Then it follows from Lemma~\ref{lem:hyp_cross_size} that this choice of $M$ satisfies
\begin{align}\label{eq:MN}
M \ge \sup_{1/ \alpha < \lambda \le 2} \left( N \prod_{j=1}^d (1 + 2 \gamma_j^\lambda \zeta(\alpha \lambda))^{-1} \right)^{1/\lambda} =: Z.
\end{align}
The function of $\lambda$ displayed in \eqref{eq:MN} extends continuously to $\lambda=1/\alpha$ with value $0$; hence this supremum, and the analogous suprema used below, are attained.

\subsection{Rank-1 lattice point set}\label{subsec:lattice}

For a prime number $N$, a rank-1 lattice point set in dimension $d$ with $N$ points is defined by a generating vector $\bsg \in \{1,2,\dots,N-1\}^d$ as
\[ 
\mathcal{L}_N(\bsg) = \left\{ \left\{ \frac{n \bsg}{N}\right\} \in [0,1)^d\mid  n = 0, 1, \ldots, N-1\right\},
\]
where $\{x \} = x - \lfloor x \rfloor$ denotes the fractional part of a real number $x\ge 0$ and is applied componentwise to a vector. For a shift vector $\bsy \in [0,1)^d$, the shifted rank-1 lattice point set is defined as
\[ 
\mathcal{L}_N(\bsg, \bsy) = \left\{ \left\{ \frac{n \bsg}{N} + \bsy\right\}  \in [0,1)^d \mid n = 0, 1, \ldots, N-1\right\} .
\]
Here, (shifted) rank-1 lattice point sets satisfy a useful property in the context of Fourier approximation.
For any $\bsk\in \ZZ^d$ and $\bsy\in [0,1)^d$, we have
\begin{align}
    \frac{1}{N}\sum_{\bsx\in \mathcal{L}_N(\bsg,\bsy)} \mathrm{e}^{2\pi \mathrm{i} \bsk \cdot \bsx} = \begin{cases}
        \mathrm{e}^{2\pi \mathrm{i} \bsk \cdot \bsy} & \text{if $\bsk\in \mathcal{L}_N^\perp(\bsg)$,}\\
        0 & \text{otherwise,}
    \end{cases}
\end{align}
where the \emph{dual lattice} $\mathcal{L}_N^\perp(\bsg)$ is defined as
\[ 
\mathcal{L}_N^\perp(\bsg) = \{ \bsk \in \mathbb{Z}^d \mid \bsg \cdot \bsk \equiv 0 \pmod{N}\}.
\]

For given $\alpha$, $\bsgamma$, and a chosen radius $M>0$, recall that the hyperbolic cross set $A_{\alpha,\bsgamma,M}$ is defined by \eqref{eq:Def_hyperbolic_set}. For any $\bsell \in \mathbb{Z}^d$,  
we define the \emph{fiber} in $A_{\alpha, \bsgamma, M}$ along $\bsell$ by
\begin{equation*}
\Gamma_{\alpha, \bsgamma, N}(\bsell; \bsg) := \{ \bsk \in A_{\alpha, \bsgamma, M}: \bsk \cdot \bsg \equiv \bsell \cdot \bsg \pmod{N} \}.
\end{equation*}
In particular, if $\bsell \in A_{\alpha, \bsgamma, M}$, then we have $\bsell \in \Gamma_{\alpha, \bsgamma, N}(\bsell; \bsg)$. Note that the fiber can also be expressed as
\begin{equation}\label{fiber2}
\Gamma_{\alpha, \bsgamma, N}(\bsell; \bsg) = A_{\alpha, \bsgamma, M}\cap \left( \bsell+\mathcal{L}_N^\perp(\bsg)\right).
\end{equation}
Furthermore, we define the fiber length $R_{\alpha, \bsgamma, N}(\bsell; \bsg)$ and the maximal fiber length $R_{\alpha,\bsgamma,N}(\bsg)$,
\[  R_{\alpha, \bsgamma, N}(\bsell; \bsg) = |\Gamma_{\alpha, \bsgamma, N}(\bsell; \bsg)|, \quad \text{and}\quad R_{\alpha, \bsgamma, N}(\bsg) = \max_{\bsell \in \mathbb{Z}^d } R_{\alpha, \bsgamma, N}(\bsell; \bsg).
\]
The maximal fiber length can also be written as $R_{\alpha, \bsgamma, N}(\bsg) = \max_{\bsell \in A_{\alpha, \bsgamma, M}} R_{\alpha, \bsgamma, N}(\bsell; \bsg)$.

Let $\bsell^{(1)}, \bsell^{(2)}, \ldots, \bsell^{(J)} \in A_{\alpha, \bsgamma, M}$ be such that
\[
\Gamma_{\alpha, \bsgamma, N}(\bsell^{(i)}; \bsg)\, \cap\, \Gamma_{\alpha, \bsgamma, N}(\bsell^{(j)}; \bsg) = \emptyset,
\]
for $i \neq j$, and 
\[
\bigcup_{j=1}^J \Gamma_{\alpha, \bsgamma, N}(\bsell^{(j)}; \bsg) = A_{\alpha, \bsgamma, M}.
\]
This partition of the hyperbolic cross into fibers is instrumental in addressing the aliasing phenomenon: each fiber corresponds to a set of frequencies that are indistinguishable under the single rank-1 lattice rule, and in our algorithm, we handle each fiber separately by applying suitable shifts to accurately estimate the corresponding Fourier coefficients.

\begin{remark}
Since the sets $\Gamma_{\alpha, \bsgamma,N}(\bsell^{(j)}; \bsg)$ partition the set $A_{\alpha, \bsgamma,M}$, we obtain
\begin{equation*}
|A_{\alpha,\bsgamma, M}| = \sum_{j=1}^J |\Gamma_{\alpha,\bsgamma, N}(\bsell^{(j)}; \bsg)| \ge R_{\alpha, \bm{\gamma}, N}(\bm{g}).
\end{equation*}
\end{remark}

\section{The approximation algorithms and the construction of shifted lattice point sets}\label{sec:algorithm}

In our approximation algorithm, we approximate the Fourier coefficients in the set $A_{\alpha, \bsgamma, M}$ using a lattice rule with generating vector $\bsg$. However, a lattice rule cannot distinguish between frequencies that lie in the same fiber $\Gamma_{\alpha,\bsgamma,N}(\bsell;\bsg)$, leading to aliasing. To address this, we employ shifted lattice rules. We begin by introducing a technique that eliminates aliasing among the Fourier coefficients.

Consider $\bsell \in A_{\alpha, \bsgamma, M}$ such that $R_{\alpha, \bsgamma, N}(\bsell; \bsg) > 1$ and let $\bsell \neq \bsk \in \Gamma_{\alpha, \bsgamma, N}(\bsell; \bsg)$. Then the frequencies $\bsk, \bsell$ are aliased since $\bsk \cdot \bsg \equiv \bsell \cdot \bsg \pmod{N}$. When approximating Fourier coefficients using this lattice rule, this means that
\begin{equation*}
\widehat{f}(\bsell) \approx \frac{1}{N}\sum_{n=0}^{N-1} f(\{ n \bsg / N\})\, \mathrm{e}^{-2\pi \mathrm{i} n \bsell \cdot \bsg / N } = \frac{1}{N}\sum_{n=0}^{N-1} f(\{ n \bsg / N\})\, \mathrm{e}^{-2\pi \mathrm{i} n \bsk \cdot \bsg / N} \approx \widehat{f}(\bsk),
\end{equation*}
i.e., all the Fourier coefficients with frequencies in the same fiber $\Gamma_{\alpha, \bsgamma,N}(\bsell; \bsg)$ are approximated by the same value, as the lattice rule cannot distinguish between them.

To avoid this problem, \cite{KV19} used multiple lattice rules; that is, if $\bm{k}_1, \bm{k}_2$ are aliased with respect to a generating vector $\bm{g}_1$, then one can use a lattice rule with a different generating vector $\bm{g}_2$ for which the frequencies $\bm{k}_1, \bm{k}_2$ are not aliased. In this paper, we use a different approach, namely we associate with each frequency $\bsell_m$ in a given fiber, say $\Gamma_{\alpha, \bsgamma, N}(\bsell; \bsg) = \{\bsell_1, \ldots, \bsell_v\}$ (where $v = R_{\alpha, \bsgamma, N}(\bsell; \bsg)$), shifts $\bsy_m^{(s)} \in [0,1)^d$, where $m = 1, 2, \ldots, v$ and $s = 1, 2, \ldots, S$. Although a square system with $S=1$ can be nonsingular for particular shifts, our analysis chooses $S>1$ in order to obtain a quantitative conditioning bound with high probability. Using approximations based on shifted lattice rules, we obtain approximations of $\widehat{f}(\bsell_m)$, $m = 1, 2, \ldots, v$, as we show in the following section.

\subsection{The deterministic approximation algorithm}\label{sec_alg_det}

We illustrate how to derive an approximation method for all the frequencies in $\bsell_1, \bsell_2, \ldots, \bsell_v \in \Gamma_{\alpha,\bsgamma,N}(\bsell; \bsg)$. 
For $m = 1, 2, \ldots, v$ and $s = 1, 2, \ldots, S$ we define\footnote{An alternative algorithm would be to use $N^{-1} \sum_{n=0}^{N-1} f(\{n \bsg/N + \bsy_m^{(s)}\}) \exp(-2\pi \mathrm{i} \bsell_m \cdot (n \bsg/N+\bsy_m^{(s)}))$. Similar results hold for this variation; however, the algorithm becomes slightly more involved.}
\begin{align*}
\mathcal{F}_{N}(f, \bsell_m, \bsy_m^{(s)}) = \frac{1}{N} \sum_{n=0}^{N-1} f(\{n \bsg / N + \bsy_m^{(s)} \})\, \mathrm{e}^{-2\pi \mathrm{i} \bsell_m \cdot ( n \bsg / N  )}   = \sum_{\bsk \in \mathcal{L}_N^\perp(\bsg)} \widehat{f}(\bsk + \bsell_m)\, \mathrm{e}^{2 \pi \mathrm{i}  (\bsk + \bsell_m)  \cdot \bsy_m^{(s)} }.
\end{align*}
Setting $\bsk = \bsell_i - \bsell_m $ and ignoring frequencies not in the fiber $\Gamma_{\alpha, \bsgamma, N}(\bsell; \bsg)$ in the above sum (to illustrate the main idea), we have  
\begin{equation*}
\mathcal{F}_N(f, \bsell_m, \bsy_m^{(s)}) \approx \sum_{i=1}^v \widehat{f}(\bsell_i)\, \mathrm{e}^{2\pi \mathrm{i}  \bsell_i \cdot \bsy_m^{(s)} }, \quad \mbox{for } m = 1, 2, \ldots, v, \mbox{ and } s = 1, 2, \ldots, S.
\end{equation*}
In matrix form, we can write this as
\begin{equation*}
\begin{pmatrix} \mathcal{F}_{N}(f,\bsell_1, \bsy_1^{(1)}) \\ 
\vdots \\ \mathcal{F}_{N}(f,\bsell_v, \bsy_v^{(S)} ) \end{pmatrix} \approx  B(\bsell, \{\bsy_m^{(s)}\}) \begin{pmatrix}  \widehat{f}(\bsell_1) \\ \vdots \\ \widehat{f}(\bsell_v) \end{pmatrix},
\end{equation*}
where we define the matrix $B(\bsell, \{\bsy_m^{(s)}\}) \in \mathbb{C}^{vS \times v}$ by
\begin{equation*}
B(\bsell, \{\bsy_m^{(s)}\}) = \left[ \mathrm{e}^{2\pi \mathrm{i} \bsell_{i} \cdot \bsy_m^{(s)}} \right]_{1 \le i, m \le v, 1 \le s \le S} =  \begin{pmatrix} \mathrm{e}^{2\pi \mathrm{i} \bsell_1  \cdot \bsy_1^{(1)}} & \ldots & \ldots  & \mathrm{e}^{2\pi \mathrm{i} \bsell_v  \cdot \bsy_1^{(1)}}  \\ \vdots & \ddots &  &   \vdots \\  \vdots & & \ddots & \vdots \\ \mathrm{e}^{2\pi \mathrm{i} \bsell_1 \cdot \bsy_v^{(S)}} & \cdots & \ldots & \mathrm{e}^{2\pi \mathrm{i} \bsell_v \cdot \bsy_v^{(S)}}  \end{pmatrix}.
\end{equation*}

\begin{remark}\rm
It is useful to distinguish two equivalent ways of viewing the shifted samples. In the definition of $\mathcal F_N(f,\boldsymbol\ell_m,\boldsymbol y_m^{(s)})$ above, the underlying lattice $\mathcal L_N(\boldsymbol g)$ is kept fixed. Thus, if we define the shifted function $f_{\boldsymbol{y}}(\boldsymbol{x}):=
f({\boldsymbol x+\boldsymbol y})$, then $\mathcal F_N(f,\boldsymbol\ell_m,\boldsymbol y_m^{(s)})$ is the discrete Fourier coefficient of $f_{\boldsymbol y_m^{(s)}}$ with respect to the fixed sampling set $\mathcal L_N(\boldsymbol g)$, namely
\[
\mathcal F_N(f,\boldsymbol\ell_m,\boldsymbol y_m^{(s)}) = 
\frac1N\sum_{n=0}^{N-1}
f_{\boldsymbol y_m^{(s)}} \left(\left\{ n\boldsymbol g/N \right\}\right) \mathrm{e}^{-2\pi \mathrm{i} \boldsymbol\ell_m\cdot (n\boldsymbol g/N)}.
\]
With this convention, the lattice and hence also the dual lattice and the fibers $\Gamma_{\alpha,\boldsymbol{\gamma},N}(\boldsymbol\ell;\boldsymbol g)$ do not depend on the shift. This leads to
\[
\mathcal F_N(f,\boldsymbol\ell_m,\boldsymbol y_m^{(s)}) = 
\sum_{\boldsymbol k\in \mathcal L_N^\perp(\boldsymbol g)}
\widehat f(\boldsymbol k+\boldsymbol\ell_m)
\mathrm{e}^{2\pi \mathrm{i} (\boldsymbol k+\boldsymbol\ell_m)\cdot \boldsymbol y_m^{(s)}},
\]
and therefore, for the aliased frequencies $\boldsymbol\ell_i \in
\Gamma_{\alpha,\boldsymbol{\gamma},N}(\boldsymbol\ell;\boldsymbol g)$, the leading terms give the matrix entries $\mathrm{e}^{2\pi \mathrm{i} \boldsymbol\ell_i\cdot \boldsymbol y_m^{(s)}}$.

Alternatively, one may regard $\mathcal L_N(\boldsymbol g)+\boldsymbol y_m^{(s)}$ as the sampling set and compute the discrete Fourier coefficient of $f$ with respect to this shifted discretisation:
\[
\widetilde{\mathcal F}_N(f,\boldsymbol\ell_m,\boldsymbol y_m^{(s)})
:=
\frac1N\sum_{n=0}^{N-1}
f\left(\left\{n\boldsymbol g / N + \boldsymbol y_m^{(s)} \right\} \right)
\mathrm{e}^{ -2\pi \mathrm{i}\boldsymbol\ell_m\cdot
( n\boldsymbol g / N + \boldsymbol y_m^{(s)} ) }.
\]
The two conventions differ only by the known phase factor
\[
\widetilde{\mathcal F}_N(f,\boldsymbol\ell_m,\boldsymbol y_m^{(s)}) =
\mathrm{e}^{ -2\pi \mathrm{i} \boldsymbol\ell_m\cdot \boldsymbol y_m^{(s)} } \mathcal F_N(f,\boldsymbol\ell_m,\boldsymbol y_m^{(s)}).
\]
Thus the same information is used in both cases, but the phase factors appear in different places in the resulting least-squares system. In this section we use the first convention because it keeps the lattice, the dual lattice, and the fibers fixed while the shifts only enter through the function values and the phase matrix.
\end{remark}

In Section~\ref{sec_good_shifts} we show that shifts $\{\bsy_m^{(s)}\}$ can be found with high probability such that the matrix $B(\bsell, \{\bsy_m^{(s)}\})^H B(\bsell, \{\bsy_m^{(s)}\})$ is invertible. Then we use the approximation method
\begin{equation}\label{eq_def_B}
\begin{pmatrix}  \widehat{f}(\bsell_1) \\ \vdots \\ \widehat{f}(\bsell_v) \end{pmatrix} \approx \left[ B(\bsell, \{\bsy_m^{(s)}\})^H B(\bsell, \{\bsy_m^{(s)}\}) \right]^{-1}  B(\bsell, \{\bsy_m^{(s)}\})^H 
\begin{pmatrix} \mathcal{F}_{N}(f,\bsell_1, \bsy_1^{(1)}) \\ \vdots \\ \mathcal{F}_{N}(f,\bsell_v, \bsy_v^{(S)} ) \end{pmatrix}
= : \begin{pmatrix}  \mathcal{B}_N(f, \bsell_1) \\ \vdots \\ \mathcal{B}_N(f, \bsell_v)  \end{pmatrix}.
\end{equation}
Finally, we define the approximation of $f$ by
\begin{equation}\label{eq_alg_det}
f \approx \mathcal{A}(f) := \sum_{j=1}^J\sum_{\bsk \in \Gamma_{\alpha,\bsgamma,N}(\bsell^{(j)},\bsg)} \mathcal{B}_N(f, \bsk)\, \mathrm{e}^{2 \pi \mathrm{i} \bsk \cdot \bsx},
\end{equation}
where $\Gamma_{\alpha,\bsgamma,N}(\bsell^{(i)}; \bsg) \cap \Gamma_{\alpha, \bsgamma,N}(\bsell^{(j)}; \bsg) = \emptyset$ for all $i \neq j$ and $\bigcup_{j=1}^J \Gamma_{\alpha,\bsgamma, N}(\bsell^{(j)}; \bsg) = A_{\alpha,\bsgamma,M}.$

The matrix $B(\bsell, \{\bsy_m^{(s)}\})^H B(\bsell, \{\bsy_m^{(s)}\})$ is of size $v \times v$ with $v \le R_{\alpha, \bsgamma, N}(\bsg)$. Therefore, it is beneficial to use generating vectors $\bsg$ for which $R_{\alpha, \bsgamma, N}(\bsg)$ is `small'. In Section~\ref{sec_exist_g} we show that we can find $\bsg$ with $R_{\alpha, \bsgamma, N}(\bsg) = \mathcal{O}( (\log N)^{d-1} )$. Note that we can use the same shifts $\{\bsy_m^{(s)}\}$ for each fiber $\Gamma_{\alpha, \bsgamma, N}(\bsell^{(j)}; \bsg)$, i.e. $\{\bsy_m^{(s)}\}$ does not depend on $j$.

The randomized algorithm is based on the same idea and is presented in the next section.

\subsection{The randomized algorithm}\label{sec_alg_rand}

Let $\bsDelta \in [0,1)^d$ be chosen uniformly at random.
Let $\bsell \in A_{\alpha, \bsgamma, M}$ and $\Gamma_{\alpha, \bsgamma, N}(\bsell; \bsg) = \{\bsell_1, \bsell_2, \ldots, \bsell_v\}$. For $m = 1, 2, \ldots, v$ we use random shifts $\{ \bsy_m^{(s)} + \bsDelta \}_{m,s}$, i.e.\footnote{An alternative algorithm would be to use $N^{-1} \sum_{n=0}^{N-1} f(\{n \bsg/N + \bsy_m^{(s)} + \bsDelta \}) \exp(-2\pi \mathrm{i} \bsell_m \cdot (n \bsg/N+ \bsDelta ))$. Similar results hold for this variation; however, in this case we get $\mathbb{E}_{\bsDelta} \mathcal{F}_N(f, \bsell_m, \bsy_m^{(s)} + \bsDelta) = \widehat{f}(\bsell_m) \mathrm{e}^{2 \pi \mathrm{i} \bm{\ell}_m \cdot \bm{y}_m^{(s)}}$, rather than $\mathbb{E}_{\bsDelta} \mathcal{F}_N(f, \bsell_m, \bsy_m^{(s)} + \bsDelta) = \widehat{f}(\bsell_m)$ as we do with the current method.}
\begin{align*}
\mathcal{F}_{N}(f, \bsell_m, \bsy_m^{(s)} + \bsDelta)  & =  \frac{1}{N} \sum_{n=0}^{N-1} f(\{n \bsg / N + \bsy_m^{(s)} + \bsDelta \})\, \mathrm{e}^{-2\pi \mathrm{i} \bsell_m \cdot ( n \bsg / N  + \bsy_m^{(s)} + \bsDelta )}  \\ &  = \sum_{\bsk \in \mathcal{L}_N^\perp(\bsg)} \widehat{f}(\bsk + \bsell_m)\, \mathrm{e}^{2 \pi \mathrm{i}  \bsk   \cdot (\bsy_m^{(s)} + \bsDelta) }.
\end{align*}

To derive an approximation method for all the frequencies in $\Gamma_{\alpha, \bsgamma, N}(\bsell; \bsg)$, by setting $\bsk = \bsell_i - \bsell_m $ in the above sum and ignoring frequencies not in the fiber $\Gamma_{\alpha,\bsgamma,N}(\bsell; \bsg)$, we have
\begin{align*}
\mathcal{F}_N(f, \bsell_m, \bsy_m^{(s)} + \bsDelta) & \approx \sum_{i=1}^v \widehat{f}(\bsell_i)\, \mathrm{e}^{2\pi \mathrm{i}  (\bsell_i - \bsell_m) \cdot ( \bsy_m^{(s)} + \bsDelta) } \\ & = \mathrm{e}^{- 2\pi \mathrm{i} \bsell_m \cdot (\bsy_m^{(s)} + \bsDelta ) } \sum_{i=1}^v \widehat{f}(\bsell_i)\, \mathrm{e}^{2 \pi \mathrm{i} \bsell_i  \cdot (\bsy_m^{(s)} +  \bsDelta )},
\end{align*}
for  $m = 1, 2, \ldots, v$, and $s = 1, 2, \ldots, S$.
Notice that $\mathbb{E}_{\bsDelta} \mathcal{F}_N(f, \bsell_m, \bsy_m^{(s)} + \bsDelta) = \widehat{f}(\bsell_m)$. Define the diagonal matrix  $$D(\bsDelta ) = (d_{i,j}(\bsDelta))_{i, j = 1, \ldots, v}, \mbox{ with } d_{i,j}(\bsDelta) = 0 \mbox{ for } i \neq j, \mbox{ and } d_{i,i}(\bsDelta) = \mathrm{e}^{2\pi \mathrm{i} \bsell_i \cdot \bsDelta} \mbox{ for } i=1, \ldots, v.$$ Then
\begin{equation*}
\begin{pmatrix} \mathrm{e}^{2\pi \mathrm{i} \bsell_1 \cdot (\bsy_1^{(1)} + \bsDelta)} \mathcal{F}_{N}(f,\bsell_1, \bsy_1^{(1)} + \bsDelta) \\ \vdots \\ \mathrm{e}^{2\pi \mathrm{i} \bsell_v \cdot (\bsy_v^{(S)} + \bsDelta ) }  \mathcal{F}_{N}(f,\bsell_v, \bsy_v^{(S)} + \bsDelta) \end{pmatrix} \approx  B(\bsell,\{\bsy_m^{(s)}\}) D(\bsDelta) \begin{pmatrix}  \widehat{f}(\bsell_1) \\ \vdots \\  \widehat{f}(\bsell_v) \end{pmatrix}.
\end{equation*}

We now use the approximation method
\begin{align*}
\begin{pmatrix}  \widehat{f}(\bsell_1) \\ \vdots \\  \widehat{f}(\bsell_v) \end{pmatrix} & \approx D(\bsDelta)^{-1} \left[ B(\bsell, \{\bsy_m^{(s)} \})^H B(\bsell, \{\bsy_m^{(s)}  \}) \right]^{-1}  B(\bsell, \{\bsy_m^{(s)}  \})^H \begin{pmatrix} \mathrm{e}^{2\pi \mathrm{i} \bsell_1 \cdot (\bsy_1^{(1)} + \bsDelta ) } \mathcal{F}_{N}(f,\bsell_1, \bsy_1^{(1)} + \bsDelta) \\ \vdots \\ \mathrm{e}^{2\pi \mathrm{i} \bsell_v \cdot (\bsy_v^{(S)} + \bsDelta )} \mathcal{F}_{N}(f,\bsell_v, \bsy_v^{(S)} + \bsDelta ) \end{pmatrix} \\ &
= : \begin{pmatrix}  \mathcal{B}_N(f, \bsell_1, \{\bsy_m^{(s)} \}, \bsDelta) \\ \vdots \\ \mathcal{B}_N(f, \bsell_v, \{\bsy_m^{(s)}\},  \bsDelta )  \end{pmatrix}.
\end{align*}

Then we define the approximation of $f$ by
\begin{equation}\label{eq_alg_rand}
f \approx \mathcal{A}_{\bsDelta}(f) := \sum_{\bsk \in A_{\alpha, \bsgamma, M}}  \mathcal{B}_N(f, \bsk, \{\bsy_m^{(s)} \}, \bsDelta )\, \mathrm{e}^{2 \pi \mathrm{i} \bsk \cdot \bsx}. 
\end{equation}
Since $\bsDelta$ is chosen randomly from the uniform distribution in $[0,1)^d$, the approximation $\mathcal{A}_{\bsDelta}$ is also random.

\subsection{Summary of the complete algorithm}

For clarity, Algorithm~\ref{alg:complete-approximation} summarizes the
complete workflow of the proposed approximation method. The deterministic
and randomized algorithms can be described simultaneously by setting
$\boldsymbol{\Delta}=\boldsymbol{0}$ in the deterministic case and choosing
$\boldsymbol{\Delta}$ uniformly from $[0,1)^d$ in the randomized case.
The formulation below is equivalent to~\eqref{eq_alg_det}
and~\eqref{eq_alg_rand}, with the phase factors incorporated
into the right-hand side of the least-squares systems.

\begin{algorithm}[H]
\caption{Shifted rank-1 lattice approximation}
\label{alg:complete-approximation}

\begin{algorithmic}[1]
\Require $f\colon[0,1)^d\to\mathbb{C}$, $\alpha>1/2$, weights $0<\gamma_j\leq1$ for $j=1,\ldots,d$, a prime number $N$, $K>1$, and the choice of deterministic or randomized approximation.
\Ensure A trigonometric approximation $\mathcal{A}_{\boldsymbol{\Delta}}(f)$.

\State Choose $M>1$ according to \eqref{eq:MN2} for the random construction or according to \eqref{eq_M} for the CBC construction, and choose a generating vector $\boldsymbol{g}\in\{1,\ldots,N-1\}^d$ as in Section~\ref{sec_exist_g} or Section~\ref{sec_g_construction}, respectively.

\State Construct the frequency index set
$
    A_{\alpha,\boldsymbol{\gamma},M}
    =
    \left\{
        \boldsymbol{k}\in\mathbb{Z}^d:
        r_{\alpha,\boldsymbol{\gamma}}(\boldsymbol{k})<M
    \right\}.
$

\State Partition $A_{\alpha,\boldsymbol{\gamma},M}$ according to
the residue $\boldsymbol{k}\cdot\boldsymbol{g}\pmod N$:
$
    A_{\alpha,\boldsymbol{\gamma},M}
    =
    \bigcup_{j=1}^{J}\Gamma_j$,
$    
    \Gamma_j
    =
    \{
        \boldsymbol{\ell}_{j,1},\ldots,
        \boldsymbol{\ell}_{j,v_j}
    \},
$
where the fibers $\Gamma_j$ are pairwise disjoint and all frequencies in
the same fiber have the same residue modulo $N$.

\State Set
$
    R=\max_{1\leq j\leq J}v_j,
    S=\left\lceil 2KR\log N \right\rceil .
$

\Repeat
    \State Draw
    $\boldsymbol{y}_m^{(s)}\sim\mathcal{U}([0,1)^d)$ independently for
    $m=1,\ldots,R$ and $s=1,\ldots,S$.

    \For{$j=1,\ldots,J$}
        \State Form the matrix $B_j\in\mathbb{C}^{v_jS\times v_j}$ with
        $
            (B_j)_{(m,s),i}
            =
            \exp\left(
                2\pi\mathrm{i}\,
                \boldsymbol{\ell}_{j,i}\cdot
                \boldsymbol{y}_m^{(s)}
            \right)$,
            $i,m=1,\ldots,v_j, 
            s=1,\ldots,S.
        $
    \EndFor
\Until{$B_j^HB_j$ is invertible and
$\|(B_j^HB_j)^{-1}\|_2\leq 1/S$ for every $j=1,\ldots,J$.}

\If{the randomized algorithm is used} 
\State Draw one additional shift
    $\boldsymbol{\Delta}\sim\mathcal{U}([0,1)^d)$.
\Else
    \State Set $\boldsymbol{\Delta}=\boldsymbol{0}$.
\EndIf

\For{$m=1,\ldots,R$ and $s=1,\ldots,S$}
    \State Evaluate $f$:
    $f_{n,m,s}
        =
        f\left(
            \left\{
                \frac{n\boldsymbol{g}}{N}
                +\boldsymbol{y}_m^{(s)}
                +\boldsymbol{\Delta}
            \right\}
        \right), n=0,\ldots,N-1.
    $
\EndFor

\For{$j=1,\ldots,J$}
    \State Form $\boldsymbol{q}_j\in\mathbb{C}^{v_jS}$ with entries
    $ (\boldsymbol{q}_j)_{(m,s)}
        =
        \frac{1}{N}\sum_{n=0}^{N-1}
        f_{n,m,s}
        \exp\left(
            -\frac{2\pi\mathrm{i}n}{N}
            \boldsymbol{\ell}_{j,m}\cdot\boldsymbol{g}
        \right),
        m=1,\ldots,v_j,
        s=1,\ldots,S.
    $

    \State Define
    $
        D_j(\boldsymbol{\Delta})
        =
        \operatorname{diag}\left(
            \exp(2\pi\mathrm{i}\,
                 \boldsymbol{\ell}_{j,1}\cdot\boldsymbol{\Delta}),
            \ldots,
            \exp(2\pi\mathrm{i}\,
                 \boldsymbol{\ell}_{j,v_j}\cdot\boldsymbol{\Delta})
        \right).
    $

    \State Recover the Fourier coefficients in the $j$-th fiber by solving
    \[
        \widetilde{\boldsymbol{f}}_j
        =
        \underset{\boldsymbol{c}\in\mathbb{C}^{v_j}}{\arg\min}\,
        \left\|
            B_j D_j(\boldsymbol{\Delta})\boldsymbol{c}
            -\boldsymbol{q}_j
        \right\|_2 .
    \]
\EndFor

\State Return
$
    \mathcal{A}_{\boldsymbol{\Delta}}(f)(\boldsymbol{x})
    =
    \sum_{j=1}^{J}\sum_{i=1}^{v_j}
    (\widetilde{\boldsymbol{f}}_j)_i
    \exp\left(
        2\pi\mathrm{i}\,
        \boldsymbol{\ell}_{j,i}\cdot\boldsymbol{x}
    \right).
$
\end{algorithmic}

\end{algorithm}

The function values associated with a fixed pair $(m,s)$ are shared by all
fibers. Moreover, the quantities $(\boldsymbol{q}_j)_{(m,s)}$ for all
residue classes can be obtained from a single length-$N$ discrete Fourier
transform of the corresponding lattice samples. In an implementation, the
least-squares problems should preferably be solved by a QR factorization
rather than by explicitly forming the inverse of $B_j^HB_j$.

\subsection{Existence of good generating vectors}\label{sec_exist_g}

In this section we show that generating vectors $\bsg$ for which $R_{\alpha,\bsgamma,N}(\bsg) = \mathcal{O}( (\log N)^{d-1} )$ can be found with high probability by randomly selecting a vector from $\{1, 2, \ldots, N-1\}^d$.

For $\bm{g} = (g_1, \ldots, g_d)\in \{1, \ldots, N-1\}^d$, we define the figure of merit
$$\rho = \rho(\bm{g} ):=\min_{\bm{\ell} \in \mathcal{L}_N^\perp(\bsg) \setminus\{\bm{0} \}} r_{\alpha, \bm{\gamma}}(\bm{\ell}) .$$

\begin{lemma}\label{lem_rho}
Let \(0<\delta<1\), $N > 2$ be a prime number and let $M \ge 1$ be such that $|A_{\alpha, \bm{\gamma}, M}| \le (N-1) (1-\delta)$. If $\bm{g}$ is sampled uniformly from \(\{1,\dots,N-1\}^d \), then
\begin{equation}\label{eq:main-bound}
\mathbb{P} \left(\rho(\bm{g}) \ge M  \right) \ge \delta.
\end{equation}
\end{lemma}

\begin{lemma}\label{lem:prob}
Let $N$ be prime and $\bm{k} \in\mathbb{Z}^s$ such that $N \nmid \bm{k}$. If $\bm{g}$ is uniform on \(\{1,\dots,N-1\}^s\), then
\[
\mathbb{P}\left( \bm{k} \cdot \bm{g} \equiv 0 \pmod N \right) \le \frac{1}{N-1} \le \frac{2}{N}.
\]
\end{lemma}

\begin{proof}
Assume that $k_j \not\equiv 0 \pmod N$. Since $N$ is prime, $k_j$ is invertible modulo $N$. Then for each choice of $g_i \in \{1, \ldots, N-1\}$, $i \in \{1, \ldots, s\} \setminus \{j\}$, there is at most one $g_j \in \{1, \ldots, N-1\}$ such that $\bm{k} \cdot \bm{g} \equiv 0 \pmod{N}$. 
\end{proof}

\begin{proof}[Proof of Lemma~\ref{lem_rho}]
Every nonzero $\bm{k}\in A_{\alpha,\bm{\gamma},M}$ has at least one component not divisible by $N$. Indeed, if $\bm{k}=N\bm{h}$ for some nonzero $\bm{h}\in\mathbb{Z}^d$, then $r_{\alpha,\bm\gamma}(q\bm h)\leq r_{\alpha,\bm\gamma}(N\bm h)<M$ for $q=1,\ldots,N$, so $\bm{h},2\bm{h},\ldots,N\bm{h}$ are $N$ distinct elements of $A_{\alpha,\bm{\gamma},M}$, contradicting $|A_{\alpha,\bm{\gamma},M}|<N$. Thus Lemma~\ref{lem:prob} applies to every summand below.
For $M \ge 1$ let
\[
X_M( \bm{g}) := \sum_{\bm{k} \in A_{\alpha, \bm{\gamma}, M} \setminus \{\bm{0}\}} \mathbf{1}\{ \bm{k} \cdot \bm{g} \equiv 0 \pmod N\}.
\]
Observe that 
\begin{equation*}
\{\bm{g}: \rho(\bm{g}) < M\} = \{\bm{g}:  X_M(\bm{g}) \ge 1\}.
\end{equation*}
Thus
\begin{equation*}
\mathbb{P}(\rho(\bm{g}) < M) \le \mathbb{E}[X_M].
\end{equation*}
Using Lemma~\ref{lem:prob}, the union bound leads to
\[
\mathbb{E}[X_M] \le \frac{|A_{\alpha, \bm{\gamma},M}| }{N-1}
 \le 1-\delta, 
\]
therefore
\begin{equation*}
\mathbb{P}( \rho(\bm{g}) \ge M) \ge \delta.
\end{equation*}
\end{proof}

\subsection{Construction of generating vectors with large figure of merit}\label{sec_g_construction}

Rather than selecting $\bm{g}$ randomly, one can also use a component-by-component construction minimizing $P_{\alpha, \bm{\gamma}, N}$ given by
\begin{equation*}
P_{\alpha, \bm{\gamma},N}(\bm{g}) = \sum_{\bm{k} \in \mathcal{L}_N^\perp(\bsg)\setminus \{\bm{0}\}} \frac{1}{r^2_{\alpha, \bm{\gamma}}(\bm{k})}.
\end{equation*}
This quantity is the squared worst-case error for the integration problem by a rank-1 lattice rule in the weighted Korobov space $\mathcal{K}_{d,\alpha,\bsgamma}$.
The largest term in this sum is $1/\rho^2(\bm{g})$, therefore
\begin{equation*}
\frac{1}{\rho^2(\bm{g})} < P_{\alpha, \bm{\gamma},N}(\bm{g}).
\end{equation*}
Since $P_{\alpha, \bm{\gamma},N}(\bm{g})$ is convergent for any $\alpha > 1/2$, we can use a component-by-component construction to find $\bm{g}$ with small $P_{\alpha, \bm{\gamma},N}(\bm{g})$, and therefore large $\rho(\bm{g})$. When $\alpha$ is a positive integer, the sum $P_{\alpha,\bsgamma,N}(\bsg)$ can be analytically computed as
\begin{align*}
P_{\alpha,\bsgamma,N}(\bsg)&=-1+\sum_{\bsk\in\mathbb{Z}^d} \frac{1}{r^2_{\alpha,\bsgamma}(\bsk)}\frac{1}{N}\sum_{n=0}^{N-1}\mathrm{e}^{2\pi i(\bsk\cdot\bsg)n/N} \\
&=-1+\frac{1}{N}\sum_{n=0}^{N-1}\prod_{j=1}^d\left( 1+\gamma_j^2 \sum_{k\in\mathbb{Z}\setminus \{0\}} \frac{1}{|k|^{2\alpha}}\mathrm{e}^{2\pi i k g_j n/N}\right)\\
&=-1+\frac{1}{N}\sum_{n=0}^{N-1}\prod_{j=1}^d\left( 1+ \gamma_j^2 \frac{(-1)^{\alpha +1}(2\pi)^{2\alpha}}{ (2\alpha )!} B_{2\alpha}\left(\left\{ \frac{n g_j}{N}  \right\} \right)\right),
\end{align*}
where $B_{2\alpha}$ denotes the Bernoulli polynomial of degree $2 \alpha$. We list the CBC construction in Algorithm \ref{cbc_principle}. Based on the work of Nuyens and Cools \cite{NC06}, computing $P_{\alpha,\bsgamma, N}((g_1,\ldots,g_s,g))$  for all $g \in \{1, \ldots,N-1\}$ (with given $g_1,\ldots,g_s$) can be done with $\mathcal{O}(N \log N)$ arithmetic operations with the help of the fast Fourier transform, and the total cost of generating $\bsg$ with a large figure of merit is $\mathcal{O}(dN\log N)$.

\begin{algorithm}[H]
\caption{CBC construction for small $P_{\alpha,\bsgamma,N}$}
\label{cbc_principle}
Let $d,\alpha\in\mathbb{N}$, $\gamma_1,\ldots,\gamma_d\in \R_{>0}$, and $N$ be a prime number. 
Construct a generating vector $\bsg=(g_1,\ldots,g_d)\in \{1,\ldots,N-1  \}^d$ as follows.
\begin{enumerate}
\item Choose $g_1=1$.
\item \textbf{For} $s$ from $1$ to $d-1$ do the following:\\
Assume that $(g_1,\ldots,g_{s})\in \{1,\ldots,N-1 \}^{s}$ have already been found. Choose $g_{s+1}\in \{ 1,\ldots,N-1 \}$ as
\begin{align*}
g_{s+1}:=\mathrm{argmin}_{g\in \{1,\ldots,N-1\}}\; P_{\alpha,\bsgamma,N}((g_1,\ldots,g_{s},g)).
\end{align*}
\textbf{end for}
\item Set $\bsg=(g_1,\ldots,g_{d}).$
\end{enumerate}
\end{algorithm}

For a generating vector $\bsg$ found by the CBC construction (Algorithm~\ref{cbc_principle}), it is well known that 
\begin{align*}
P_{\alpha,\bsgamma,N}(\bsg)\leq \left( \frac{2}{N} \prod_{j=1}^d \left(1+2 \gamma_j^{2\lambda'} \zeta(2\alpha \lambda' )\right) \right)^{1/\lambda'},
\end{align*}
for any $\lambda' \in (1/(2\alpha), 1]$ (cf.\ \cite[Theorem~3.7]{DKP22} in which $\gamma^2_j$ is replaced by $\gamma_j$). With $\lambda = 2\lambda'$, we have
\[ P_{\alpha,\bsgamma,N}(\bsg)\leq \left( \frac{2}{N} \prod_{j=1}^d \left(1+2 \gamma_j^{\lambda} \zeta(\alpha \lambda)\right) \right)^{2/\lambda}, \]
for any $\lambda \in (1/\alpha,2]$. Since $1/(\rho(\bsg))^2< P_{\alpha,\bsgamma,N}(\bsg)$, we obtain
\begin{equation}\label{cbc_rho}
\rho(\bsg)> \max_{1/\alpha < \lambda \le 2} \left( \frac{N}{2} \prod_{j=1}^d(1+2\gamma_j^{\lambda} \zeta(\alpha \lambda))^{-1} \right) ^{1/\lambda}.
\end{equation}

\subsection{A bound on the fiber length}

With the bound on $\rho(\bm g)$, we now prove directly a bound
on the maximal fiber length $R_{\alpha,\bm\gamma,N}(\bm g)$.

\begin{lemma}\label{lem_g_exist}
Let $N>2$ be prime, let $\delta\in(0,1)$, assume that $0<\gamma_j\leq1$ for $j=1,\ldots,d$, and let $M>1$ satisfy $|A_{\alpha,\bm\gamma,M}|\leq(N-1)(1-\delta)$. Let $m\in\mathbb N$ satisfy $2^{m-1}<M^{1/\alpha}\leq2^m$, and choose $\bm g$ uniformly from $\{1,\ldots,N-1\}^d$. Then, with probability at least $\delta$, $\rho(\bm g)\geq M$ and
\begin{equation}\label{bound_g_exist}
R_{\alpha,\bm\gamma,N}(\bm g)
\leq 1+(3^d-1)(m+d)^{d-1}.
\end{equation}

If $\alpha\in\mathbb N$, $\bm g$ is constructed by Algorithm~\ref{cbc_principle}, the value of $M$ below satisfies $M>1$, and
\begin{equation}\label{eq_M}
M  = \max_{1/\alpha < \lambda \le 2} \left(\frac{N}{2} \prod_{j=1}^d (1 + 2 \gamma_j^\lambda \zeta(\alpha \lambda))^{-1} \right)^{1/\lambda},
\end{equation}
then $|A_{\alpha,\bm\gamma,M}|\leq N/2$, $\rho(\bm g)>M$, and \eqref{bound_g_exist} holds for the integer $m$ determined by $2^{m-1}<M^{1/\alpha}\leq2^m$.
\end{lemma}

\begin{proof}
By Lemma~\ref{lem_rho}, the event $\rho(\bm g)\geq M$ has
probability at least $\delta$.  We prove the fiber bound on this event.
Fix $\bm\ell\in\mathbb Z^d$ and put
$S_{\bm\ell}=\bm\ell+\mathcal L_N^\perp(\bm g)$.  By the definition of
$\rho(\bm g)$, distinct points $\bm k,\bm k'\in S_{\bm\ell}$ satisfy
\begin{equation}\label{eq_r_rho_M}
r_{\alpha,\bm\gamma}(\bm k-\bm k')\geq\rho(\bm g)\geq M.
\end{equation}

The zero vector contributes at most one point.  Fix a nonempty
support $u\subseteq\{1,\ldots,d\}$, put $q=|u|$, fix the signs of the
coordinates in $u$, and choose one index $j_0\in u$.  For
$\bm k\in A_{\alpha,\bm\gamma,M}$ having this support and these signs,
set
\[
t_j(\bm k)=\frac{|k_j|}{\gamma_j^{1/\alpha}}\qquad(j\in u).
\]
Since $0<\gamma_j\leq1$ and $k_j\neq0$, we have $t_j(\bm k)\geq1$.
For every $j\in u\setminus\{j_0\}$, let $n_j\geq1$ be determined by
\[
2^{n_j-1}\leq t_j(\bm k)<2^{n_j}.
\]
The inequality $r_{\alpha,\bm\gamma}(\bm k)<M$ implies
\[
2^{\sum_{j\in u\setminus\{j_0\}}(n_j-1)}
\leq\prod_{j\in u\setminus\{j_0\}}t_j(\bm k)
<M^{1/\alpha}\leq2^m.
\]
Hence
$\sum_{j\in u\setminus\{j_0\}}(n_j-1)\leq m-1$, so the number of
possible dyadic index vectors $(n_j)_{j\in u\setminus\{j_0\}}$ is
\[
\binom{m+q-2}{q-1},
\]
with the usual interpretation that this number is $1$ when $q=1$.

For each such dyadic index vector there is at most one point of
$S_{\bm\ell}$.  Indeed, suppose that $\bm k$ and $\bm k'$ were two
distinct such points.  For $j\in u\setminus\{j_0\}$, their common sign
and common dyadic interval give
\[
\frac{|k_j-k'_j|}{\gamma_j^{1/\alpha}}<2^{n_j-1}.
\]
Choose $\bm h\in\{\bm k,\bm k'\}$ so that
$t_{j_0}(\bm h)=\max\{t_{j_0}(\bm k),t_{j_0}(\bm k')\}$.  Since the
$j_0$-coordinates also have the same sign,
\[
\frac{|k_{j_0}-k'_{j_0}|}{\gamma_{j_0}^{1/\alpha}}
<t_{j_0}(\bm h)
\]
whenever this difference is nonzero.  Omitting zero factors and using
$2^{n_j-1}\leq t_j(\bm h)$ therefore yields
\[
r_{\alpha,\bm\gamma}(\bm k-\bm k')^{1/\alpha}
<t_{j_0}(\bm h)
  \prod_{j\in u\setminus\{j_0\}}2^{n_j-1}
\leq\prod_{j\in u}t_j(\bm h)
=r_{\alpha,\bm\gamma}(\bm h)^{1/\alpha}
<M^{1/\alpha},
\]
contradicting \eqref{eq_r_rho_M}. Thus leaving one coordinate unbinned gives the
required separation.

There are $\binom dq$ choices of the support and $2^q$ choices
of its signs.  Consequently,
\[
|S_{\bm\ell}\cap A_{\alpha,\bm\gamma,M}|
\leq1+\sum_{q=1}^d\binom dq2^q\binom{m+q-2}{q-1}
\leq1+(3^d-1)(m+d)^{d-1},
\]
because $\binom{m+q-2}{q-1}\leq(m+d)^{d-1}$ and
$\sum_{q=1}^d\binom dq2^q=3^d-1$.  Taking the maximum over
$\bm\ell$ proves \eqref{bound_g_exist}.

In the CBC case, let $\lambda$ be a maximizer in \eqref{eq_M}.
Lemma~\ref{lem:hyp_cross_size} then gives
\[
|A_{\alpha,\bm\gamma,M}|
\leq M^\lambda\prod_{j=1}^d
       \bigl(1+2\gamma_j^\lambda\zeta(\alpha\lambda)\bigr)
=\frac N2,
\]
while \eqref{cbc_rho} gives $\rho(\bm g)>M$.  Thus the preceding fiber
bound applies and completes the proof.
\end{proof}

For either choice of $M$ in \eqref{eq_M} or \eqref{eq:MN2}, the product appearing there is at least $1$ and $1/\lambda<\alpha$. Hence $M\leq N^\alpha$, so $M^{1/\alpha}\leq N$ and $m<1+\log_2N$. In particular, $m=\mathcal O(\log N)$, and Lemma~\ref{lem_g_exist} yields
\[
R_{\alpha,\bm\gamma,N}(\bm g)=\mathcal O((\log N)^{d-1}).
\]

\subsection{Existence of good shifts}\label{sec_good_shifts}

In the following, assume that $\bsg$ is a generating vector satisfying $R_{\alpha,\bsgamma,N}(\bsg) \le C (\log N)^{d-1}$, whose existence was shown in Lemma~\ref{lem_g_exist}.  Let $\bsell^{(1)}, \bsell^{(2)}, \ldots, \bsell^{(J)}$ be such that $\Gamma_{\alpha,\bsgamma,N}(\bsell^{(i)}; \bsg) \cap \Gamma_{\alpha, \bsgamma,N}(\bsell^{(j)}; \bsg) = \emptyset$ for all $i \neq j$ and $\bigcup_{j=1}^J \Gamma_{\alpha,\bsgamma, N}(\bsell^{(j)}; \bsg) = A_{\alpha,\bsgamma,M},$ and $v_j = |\Gamma_{\alpha, \bsgamma, N}(\bsell^{(j)}; \bsg)|$.  

Consider a fiber $\Gamma_{\alpha,\bsgamma,N}(\bsell^{(j)}; \bsg) = \{\bsell_1, \ldots, \bsell_{v_j}\}$.
To prove the existence of shifts $\{\bsy_m^{(s)}\}$ such that $\left\| \left[B(\bsell^{(j)}, \{\bsy_m^{(s)}\})^H B(\bsell^{(j)}, \{\bsy_m^{(s)}\})\right]^{-1} \right \|_2 \le 1/S$,  we need a lower bound on the smallest eigenvalue of the matrix
\begin{equation*}
G := B(\bsell^{(j)}, \{\bsy_m^{(s)}\})^H B(\bsell^{(j)}, \{\bsy_m^{(s)}\}) = \left[ \sum_{m=1}^{v_j} \sum_{s=1}^S \mathrm{e}^{-2\pi \mathrm{i} (\bsell_i- \bsell_{i'}) \cdot \bsy_m^{(s)} } \right]_{i, i' = 1, 2, \ldots, v_j}.
\end{equation*}
Define the entry in row $i$ and column $i'$ by
\begin{equation}\label{def_g}
g(\bsell_i, \bsell_{i'}) = \sum_{m=1}^{v_j} \sum_{s=1}^S \mathrm{e}^{-2\pi \mathrm{i} (\bsell_i- \bsell_{i'}) \cdot \bsy_m^{(s)} }.
\end{equation}
The diagonal elements of $G$ are given by $g(\bsell_i, \bsell_i) = \sum_{m=1}^{v_j} \sum_{s=1}^S \mathrm{e}^{- 2\pi \mathrm{i} (\bsell_i - \bsell_i) \cdot \bsy_{m}^{(s)}} = v_j S$.

Gershgorin's circle theorem tells us that every eigenvalue $\lambda$ lies in one of the disks
\begin{equation*}
|\lambda - v_j S| = |\lambda - g(\bsell_i, \bsell_{i}) | \le \sum_{\substack{i'=1 \\ i' \neq i}}^{v_j} | g(\bsell_i, \bsell_{i'}) |   \le (v_j -1) \max_{i \neq i'} \left| g(\bsell_i, \bsell_{i'}) \right|,
\end{equation*}
hence
\begin{equation}\label{bound_eig}
\lambda \ge v_j S - (v_j -1) \max_{i \neq i'} |g(\bsell_i, \bsell_{i'})|.    
\end{equation}

\begin{lemma}\label{lem_eig}
Assume that $|A_{\alpha,\bm\gamma,M}|<N$, let $R=R_{\alpha,\bm\gamma,N}(\bm g)$, let $K>1$, and set
\[
S=\left\lceil2KR\log N\right\rceil.
\]
If $2N^{1-K}R^2<1$, then there exist shifts $\bm y_m^{(s)}\in[0,1)^d$, $m=1,\ldots,R$, $s=1,\ldots,S$, such that, for every fiber $\Gamma_{\alpha,\bm\gamma,N}(\bm\ell^{(j)};\bm g)=\{\bm\ell_1,\ldots,\bm\ell_{v_j}\}$ and all $i\neq i'$,
\[
\left|g(\bm\ell_i,\bm\ell_{i'})\right|
=\left|\sum_{m=1}^{v_j}\sum_{s=1}^S
\mathrm e^{-2\pi\mathrm i(\bm\ell_i-\bm\ell_{i'})\cdot\bm y_m^{(s)}}\right|\leq S.
\]
If the shifts are chosen independently and uniformly from $[0,1)^d$, this holds simultaneously for all fibers with probability at least
\[
1-2N^{1-K}R^2.
\]
\end{lemma}

\begin{proof}
Fix a fiber and an ordered pair $i\neq i'$, and set $\bm h=\bm\ell_i-\bm\ell_{i'}\neq\bm0$. Then
\[
X_{\bm h}=\sum_{m=1}^{v_j}\sum_{s=1}^S
\mathrm e^{-2\pi\mathrm i\bm h\cdot\bm y_m^{(s)}}
\] 
is a sum of $v_jS$ independent, mean-zero random variables of modulus one. Identifying $\mathbb C$ with the real Hilbert space $\mathbb R^2$ and applying the Hilbert-space Hoeffding inequality to the associated partial-sum martingale (see  \cite[Theorem~3.5]{Pinelis1994}), we obtain, for $t>0$,
\[
\mathbb P(|X_{\bm h}|\geq t)
\leq2\exp\left(-\frac{t^2}{2v_jS}\right).
\]
With $t=S$ and $v_j\leq R$, this is at most $2\exp(-S/(2R))$. The number of ordered pairs over all fibers is at most $|A_{\alpha,\bm\gamma,M}|R^2<NR^2$. Hence a union bound gives a failure probability at most
\[
2NR^2\exp\left(-\frac{S}{2R}\right)
\leq2N^{1-K}R^2.
\]
Since $R=\mathcal O((\log N)^{d-1})$ in the setting of this section and $K>1$, the last quantity is smaller than $1$ for all sufficiently large $N$. Therefore suitable shifts exist, and the stated probability estimate follows.
\end{proof}

Lemma~\ref{lem_eig} shows the existence of points $\{\bsy_m^{(s)}\}$ for which $|g(\bsell_i, \bsell_{i'})|$ is always bounded by $S$. In particular, \eqref{bound_eig} implies that every eigenvalue $\lambda$ of $G$ is at least $v_j S - (v_j-1)S   = S  \ge 1$. In turn we obtain the existence of a point set $\{\bsy_m^{(s)}\}$ such that
\begin{equation*}
\left\| \left[B(\bsell^{(j)}, \{\bsy_m^{(s)}\})^H B(\bsell^{(j)}, \{\bsy_m^{(s)}\})\right]^{-1} \right \|_2 \le \frac{1}{S}  \le 1.
\end{equation*}

\begin{lemma}\label{bound_inv_mat}
Assume that $|A_{\alpha,\bm\gamma,M}|<N$ and let $S=\lceil2K R_{\alpha,\bm\gamma,N}(\bm g)\log N\rceil$ for some $K>1$. If $2N^{1-K}R_{\alpha,\bm\gamma,N}(\bm g)^2<1$, then there exist shifts such that, for every fiber,
\begin{equation*}
\left\| \left[B(\bsell^{(j)}, \{\bsy_m^{(s)}\})^H B(\bsell^{(j)}, \{\bsy_m^{(s)}\})\right]^{-1} \right\|_2 \le \frac{1}{S} \le 1.
\end{equation*}
For independently and uniformly chosen shifts, this bound holds simultaneously for all fibers with probability at least $1-2N^{1-K}R_{\alpha,\bm\gamma,N}(\bm g)^2$.
\end{lemma}

\begin{remark}
The total number $p$ of points used in the algorithms $\mathcal{A}$ and $\mathcal{A}_{\bsDelta}$ satisfies
\begin{equation}\label{no_points}
N < p = N R_{\alpha,\bsgamma,N}(\bsg) S \le N R_{\alpha, \bm{\gamma}, N}(\bm{g}) \left\lceil2 K R_{\alpha, \bm{\gamma}, N}(\bm{g}) \log N\right\rceil  \le C K N (\log N)^{2d-1}
\end{equation}
for some constant $C > 0$ independent of $N, p$, and therefore
\begin{equation}\label{bound_Np}
\frac{p}{CK (\log p)^{2d-1}} \le N < p.
\end{equation}
This controls the sampling overhead $p/N$ relative to the underlying lattice size. It does not by itself bound the oversampling factor $p/|A_{\alpha,\bsgamma,M}|$ unless $|A_{\alpha,\bsgamma,M}|$ is also known to be comparable to $N$.
Thus the consequences of Lemma~\ref{lem_g_exist} used below are, for fixed $d$,
\[
R=\mathcal O((\log N)^{d-1}),\qquad
S=\mathcal O((\log N)^d),\qquad
RS=\mathcal O((\log N)^{2d-1}),
\]
where $RS$ is the number of shifted lattice copies shared by all fibers and $p=NRS$.
\end{remark}

\section{Error analysis}\label{sec:error_analysis}

In the following two sections, we prove bounds on the approximation error in the deterministic $L_\infty$-norm setting and in the randomized $L_2$-norm setting.

\subsection{Convergence analysis for the worst-case error in $L_\infty$ norm}

We would like to analyze the worst-case $L_\infty$-approximation error as
\begin{equation*}
\mathrm{e}_{\infty}^{\wor}(\mathcal{A},\mathcal{K}_{d,\alpha,\bsgamma}) = \sup_{\substack{f \in \mathcal{K}_{d,\alpha, \bsgamma} \\ \|f\|_{d,\alpha, \bsgamma} \le 1}}\sup_{\bsx\in[0,1)^d} |f(\bsx) - \mathcal{A}(f)(\bsx)|.
\end{equation*}

We need the following lemma, whose proof can be found, for instance, in \cite[Lemma~2.5]{CG26}. In the following we write $B_j := B(\bsell^{(j)}, \{ \bsy_m^{(s)}\})$.

\begin{lemma}\label{lem_sum}
Let $\alpha>1/2$ and $M\geq 1$ be given. For any $\lambda \in (1/\alpha,2)$ we have 
\begin{align*}
\sum_{\bsk\in\mathbb{Z}^d\backslash A_{\alpha,\bsgamma,M}}\frac{1}{r_{\alpha,\bsgamma}^2(\bsk)}\leq \frac{1}{M^{2-\lambda}}\cdot\frac{8(3-\lambda)}{2-\lambda}\prod_{j=1}^d\left(1+2\gamma_j^{\lambda}\zeta(\alpha \lambda)\right) .
\end{align*}
\end{lemma}

\begin{theorem}\label{thm_det}
Let $\alpha>1/2$, let $0<\gamma_j\leq1$ for $j=1,\ldots,d$, let $N>2$ be prime, and let $\delta\in(0,1)$ and $K>1$. Choose $M$ as in \eqref{eq:MN2}, assume that $M>1$, let $m\in\mathbb N$ satisfy $2^{m-1}<M^{1/\alpha}\leq2^m$, and set
\[
R_*=1+(3^d-1)(m+d)^{d-1}.
\]
Choose $\bm g$ uniformly from $\{1,\ldots,N-1\}^d$. Conditional on $\bm g$, choose the shifts independently and uniformly from $[0,1)^d$, with
\[
S=\left\lceil2K R_{\alpha,\bm\gamma,N}(\bm g)\log N\right\rceil.
\]
By Lemma~\ref{lem_g_exist}, with probability at least $\delta$ the generating vector satisfies $\rho(\bm g)\geq M$ and $R_{\alpha,\bm\gamma,N}(\bm g)\leq R_*$. Since $m=\mathcal O(\log N)$, we have $2N^{1-K}R_*^2<1$ for all sufficiently large $N$. For such $N$, conditional on a generating vector with these properties, Lemma~\ref{bound_inv_mat} shows that the shifts satisfy $\|[B_j^HB_j]^{-1}\|_2\leq1/S$ for all fibers with probability at least $1-2N^{1-K}R_*^2$. Consequently, the two conditions hold jointly with probability at least
\[
\delta\bigl(1-2N^{1-K}R_*^2\bigr),
\]
which is positive. On this joint event, for every $\tau\in(0,1]$,
\begin{equation}\label{bound_Linf}
\mathrm{e}_{\infty}^{\wor}(\mathcal{A},\mathcal{K}_{d,\alpha,\bsgamma}) \leq  C_{d,\alpha,\bm\gamma,\delta,K,\tau}\,p^{-\alpha + 1/2 + \tau},
\end{equation}
where $\mathcal{A}(f)$ is given by \eqref{eq_alg_det} and $p$ denotes the total number of function evaluations. 

If, in addition, $\alpha\in\mathbb N$, $\bm g$ is constructed using Algorithm~\ref{cbc_principle}, $M$ is chosen as in \eqref{eq_M} and satisfies $M>1$, and the shifts satisfy $\|[B_j^HB_j]^{-1}\|_2\leq1/S$, then \eqref{bound_Linf} also holds.
\end{theorem}

\begin{proof}
The error can be separated into two parts
\begin{align}\label{wce_split}
\|f - \mathcal{A}(f) \|_{L_{\infty}}& \leq \sum_{\bsk \in A_{\alpha, \bsgamma,M}} |\widehat{f}(\bsk) - \mathcal{B}_N(f, \bsk)| + \sum_{\bsk \in \mathbb{Z}^d \setminus A_{\alpha, \bsgamma,M}} |\widehat{f}(\bsk)|.
\end{align}

For the truncation error $\sum_{\bsk \in \mathbb{Z}^d \setminus A_{\alpha, \bsgamma,M}} |\widehat{f}(\bsk)|$, we have 
\begin{align}\label{eq:truncation_error}
\sum_{\bsk \in \mathbb{Z}^d \setminus A_{\alpha, \bsgamma,M}} |\widehat{f}(\bsk)| & \le \left(\sum_{\bsk \in \mathbb{Z}^d \setminus A_{\alpha,\bsgamma, M}} |\widehat{f}(\bsk)|^2(r_{\alpha,\bsgamma}(\bsk))^2\right)^{1/2}\left(\sum_{\bsk \in \mathbb{Z}^d \setminus A_{\alpha, \bsgamma,M}} \frac{1}{(r_{\alpha,\bsgamma}(\bsk))^2}\right)^{1/2} \nonumber\\
&\leq \Vert f\Vert_{d,\alpha,\bsgamma}\left(\sum_{\bsk \in \mathbb{Z}^d \setminus A_{\alpha, \bsgamma,M}} \frac{1}{(r_{\alpha,\bsgamma}(\bsk))^2}\right)^{1/2}.
\end{align}
The last factor on the right side of \eqref{eq:truncation_error} can be estimated by Lemma~\ref{lem_sum}. 

For the aliasing error $\sum_{\bsk \in A_{\alpha,\bsgamma, M}} |\widehat{f}(\bsk) - \mathcal{B}_N(f, \bsk)|,$  we have 
\begin{align}\label{eq:aliasing_error}
\sum_{\bsk \in A_{\alpha,\bsgamma, M}} |\widehat{f}(\bsk) - \mathcal{B}_N(f, \bsk)|  & = \sum_{j=1}^J \sum_{\bsk \in \Gamma_{\alpha,\bsgamma, N}(\bsell^{(j)}; \bsg)} |\widehat{f}(\bsk) - \mathcal{B}_N(f, \bsk)|   \\ & = \sum_{j=1}^J \left\|  \begin{pmatrix}  \widehat{f}(\bsell_1) \\ \vdots \\ \widehat{f}(\bsell_{v_j}) \end{pmatrix} - \left[B_j^H B_j \right]^{-1} B_j^{H}  \begin{pmatrix} \mathcal{F}_{N}(f,\bsell_1, \bsy_1^{(1)} ) \\ \vdots \\ \mathcal{F}_{N}(f,\bsell_{v_j}, \bsy_{v_j}^{(S)}) \end{pmatrix} \right\|_1,
\end{align}
where
\begin{align*}\label{def_F_alg}
\mathcal{F}_{N}(f, \bsell_m, \bsy_m^{(s)} )  & =  \frac{1}{N} \sum_{n=0}^{N-1} f(\{n \bsg / N + \bsy_m^{(s)} \}) \mathrm{e}^{-2\pi \mathrm{i} \bsell_m \cdot  n \bsg / N }  = \sum_{\bsk - \bsell_m \in \mathcal{L}_N^\perp(\bsg)} \widehat{f}(\bsk) \mathrm{e}^{2\pi \mathrm{i} \bsk \cdot \bsy_m^{(s)}}.
\end{align*}

We have
\begin{align*}
\lefteqn{ \left[ B_j^H B_j \right]^{-1}  B_j^H \left[ \begin{pmatrix} \mathcal{F}_{N}(f,\bsell_1, \bsy_1^{(1)} ) \\ \vdots \\ \mathcal{F}_{N}(f,\bsell_{v_j}, \bsy_{v_j}^{(S)}) \end{pmatrix} - B_j \begin{pmatrix}  \widehat{f}(\bsell_1) \\ \vdots \\ \widehat{f}(\bsell_{v_j}) \end{pmatrix} \right] } \\  & = \left[ B_j^H B_j \right]^{-1}  B_j^H  \left( \sum_{\substack{ \bsk - \bsell_m \in \mathcal{L}_N^\perp(\bsg)  \\ \bsk \notin \Gamma_{\alpha, \bsgamma,N}(\bsell_m; \bsg) } } \widehat{f}(\bsk) \mathrm{e}^{2\pi \mathrm{i} \bsk  \cdot \bsy_m^{(s)} } \right)_{m = 1, \ldots, v_j, s = 1, \ldots, S} \\ & = \left[ B_j^H B_j \right]^{-1}  B_j^H \left( \sum_{\substack{ \bsk - \bsell^{(j)} \in \mathcal{L}_N^\perp(\bsg)  \\ \bsk \notin \Gamma_{\alpha,\bsgamma, N}(\bsell^{(j)}; \bsg) } } \widehat{f}(\bsk) \mathrm{e}^{2\pi \mathrm{i} \bsk  \cdot \bsy_m^{(s)} } \right)_{m = 1, \ldots, v_j, s= 1, \ldots, S} \\ & =  \sum_{\substack{ \bsk - \bsell^{(j)} \in \mathcal{L}_N^\perp(\bsg)  \\ \bsk \notin \Gamma_{\alpha,\bsgamma, N}(\bsell^{(j)}; \bsg) } } \widehat{f}(\bsk) \left[ B_j^H B_j \right]^{-1}  B_j^H \underbrace{\begin{pmatrix} \mathrm{e}^{2\pi \mathrm{i} \bsk  \cdot \bsy_1^{(1)} } \\ \vdots \\ \mathrm{e}^{2\pi \mathrm{i} \bsk \cdot \bsy_{v_j}^{(S)} } \end{pmatrix} }_{=: \bsa_{\bsk}}.
\end{align*}
For $Z_j:= \{ \bsk \in \mathbb{Z}^d: \bsk - \bsell^{(j)} \in \mathcal{L}^\perp_N(\bsg)\},$ let
\begin{equation*}
e_j(\bsk) = \left[ B_j^H B_j \right]^{-1}  B_j^H \bsa_{\bsk},
\end{equation*}
and
$E_j$ be the $v_j \times \infty$ matrix
\begin{equation*}
E_j = \left[ \frac{e_j(\bsk)}{r_{\alpha,\bsgamma}(\bsk)} \right]_{\bsk \in Z_j \setminus A_{\alpha, \bsgamma,M}}
\end{equation*}
and
\begin{equation*}
F_j = \left[ \widehat{f}(\bsk) r_{\alpha,\bsgamma}(\bsk) \right]_{\bsk \in Z_j \setminus A_{\alpha,\bsgamma, M}}.
\end{equation*}
Thus,
\begin{align*}
\sum_{\bsk \in A_{\alpha,\bsgamma, M}} |\widehat{f}(\bsk) - \mathcal{B}_N(f, \bsk)|
&=\sum_{j=1}^J\|E_jF_j\|_1 \\
&\leq \sum_{j=1}^J\sqrt{v_j}\,\|E_jF_j\|_2\\
&\leq\left(\sum_{j=1}^Jv_j\operatorname{trace}(E_j^HE_j)\right)^{1/2}
\left(\sum_{j=1}^J\|F_j\|_2^2\right)^{1/2}.
\end{align*}
Here
\begin{equation*}
\operatorname{trace}(E_j^H E_j) = \sum_{\bsk \in Z_j \setminus A_{\alpha,\bsgamma, M}} \frac{\|e_j(\bsk)\|_2^2}{r^2_{\alpha,\bsgamma}(\bsk)}.
\end{equation*}
Lemma~\ref{bound_inv_mat} implies that there are shifts such that $\|[B_j^HB_j]^{-1}\|_2\leq1/S$. Since every component of $B_j^H\bsa_{\bsk}$ has modulus at most $v_jS$,
\[
\|e_j(\bsk)\|_2
\leq \frac1S\,S v_j^{3/2}
=v_j^{3/2}
\leq R_{\alpha,\bsgamma,N}(\bsg)^{3/2}.
\]
The sets $Z_j$ are disjoint, so $\sum_{j=1}^J\|F_j\|_2^2\leq\|f\|_{d,\alpha,\bsgamma}^2$, and
\begin{equation*}
\sum_{j=1}^Jv_j\operatorname{trace}(E_j^HE_j)
\leq R_{\alpha,\bsgamma,N}(\bsg)^4
\sum_{\bsk\in\mathbb Z^d\setminus A_{\alpha,\bsgamma,M}}
\frac1{r_{\alpha,\bsgamma}(\bsk)^2}.
\end{equation*}
Consequently,
\begin{equation}\label{wce_est_main}
\sum_{\bsk \in A_{\alpha,\bsgamma, M}} |\widehat{f}(\bsk) - \mathcal{B}_N(f, \bsk)|
\leq \|f\|_{d,\alpha,\bsgamma}\,
R_{\alpha,\bsgamma,N}(\bsg)^2
\left(\sum_{\bsk\in\mathbb Z^d\setminus A_{\alpha,\bsgamma,M}}
\frac1{r_{\alpha,\bsgamma}(\bsk)^2}\right)^{1/2}.
\end{equation}

Using \eqref{wce_split}, \eqref{eq:truncation_error}, and \eqref{wce_est_main} and Lemma~\ref{lem_sum} we obtain
\begin{align*}
\|f - \mathcal{A}(f) \|_{L_{\infty}} & \le \|f\|_{d, \alpha, \bsgamma} \left[1 + R_{\alpha,\bsgamma,N}(\bsg)^2\right] \left( \sum_{\bsk \in \mathbb{Z}^d \setminus A_{\alpha,\bsgamma, M}} \frac{1}{r^2_{\alpha,\bsgamma}(\bsk)} \right)^{1/2}. 
\end{align*}
By choosing
\begin{equation}\label{eq:MN2}
M = \max_{1/\alpha < \lambda \le 2} \left( (N-1) (1-\delta) \prod_{j=1}^d (1 + 2 \gamma_j^\lambda \zeta(\alpha \lambda))^{-1} \right)^{1/\lambda}
\end{equation}
we obtain $|A_{\alpha, \bm{\gamma}, M}| \le (N-1) (1-\delta)$. Fix $\lambda\in(1/\alpha,2)$. From \eqref{eq:MN2}, $M\geq c_{\lambda,d,\alpha,\bm\gamma,\delta}N^{1/\lambda}$, and Lemma~\ref{lem_sum} together with \eqref{bound_g_exist} gives
\[
\mathrm e_\infty^{\wor}(\mathcal A,\mathcal K_{d,\alpha,\bm\gamma})
\leq C_{\lambda,d,\alpha,\bm\gamma,\delta}
(\log N)^{2d-2}N^{-1/\lambda+1/2}.
\]
Using \eqref{bound_Np}, this is bounded by $C p^{-1/\lambda+1/2}(\log p)^a$ for a constant $a>0$. Given $\tau>0$, choose $\lambda$ sufficiently close to $1/\alpha$ that $1/\lambda>\alpha-\tau/2$, and absorb $(\log p)^a$ into $p^{\tau/2}$. This proves \eqref{bound_Linf}.

The result for $\bm{g}$ constructed by Algorithm~\ref{cbc_principle} follows similarly.

\end{proof}

\subsection{Convergence analysis for the random case error}

Since the algorithm is randomized, we consider the random case error
\begin{equation*}
\mathrm{e}_2^{\ran}((\mathcal{A}_{\bsDelta})_{\bsDelta},\mathcal{K}_{d,\alpha, \bsgamma}) = \sup_{\substack{f \in \mathcal{K}_{d,\alpha, \bsgamma} \\ \|f\|_{d,\alpha, \bsgamma} \le 1}} \sqrt{ \mathbb{E}_{\bsDelta} \|f - \mathcal{A}_{\bsDelta}(f)\|_2^2 }.
\end{equation*}

\begin{theorem}\label{thm:L_2_rce}
Let $\alpha>1/2$, let $0<\gamma_j\leq1$ for $j=1,\ldots,d$, let $N>2$ be prime, and let $\delta\in(0,1)$ and $K>1$. Choose $M$ as in \eqref{eq:MN2}, assume that $M>1$, let $m\in\mathbb N$ satisfy $2^{m-1}<M^{1/\alpha}\leq2^m$, and set
\[
R_*=1+(3^d-1)(m+d)^{d-1}.
\]
Choose $\bm g$ uniformly from $\{1,\ldots,N-1\}^d$. Conditional on $\bm g$, choose the shifts independently and uniformly from $[0,1)^d$, with
\[
S=\left\lceil2K R_{\alpha,\bm\gamma,N}(\bm g)\log N\right\rceil.
\]
By Lemma~\ref{lem_g_exist}, with probability at least $\delta$ the generating vector satisfies $\rho(\bm g)\geq M$ and $R_{\alpha,\bm\gamma,N}(\bm g)\leq R_*$. Since $m=\mathcal O(\log N)$, we have $2N^{1-K}R_*^2<1$ for all sufficiently large $N$. For such $N$, conditional on a generating vector with these properties, Lemma~\ref{bound_inv_mat} shows that the shifts satisfy $\|[B_j^HB_j]^{-1}\|_2\leq1/S$ for all fibers with probability at least $1-2N^{1-K}R_*^2$. Consequently, the two conditions hold jointly with probability at least
\[
\delta\bigl(1-2N^{1-K}R_*^2\bigr),
\]
which is positive. On this joint event, for every $\tau\in(0,1]$,
\begin{equation}\label{bound_rand}
\mathrm{e}_2^{\ran}((\mathcal{A}_{\bsDelta})_{\bsDelta},\mathcal{K}_{d,\alpha, \bsgamma}) \le  C_{d,\alpha,\bm\gamma,\delta,K,\tau}\,p^{-\alpha + \tau},
\end{equation}
where $\mathcal{A}_{\bsDelta}$ is given by \eqref{eq_alg_rand} and $p$ denotes the total number of function evaluations.

If, in addition, $\alpha\in\mathbb N$, $\bm g$ is constructed using Algorithm~\ref{cbc_principle}, $M$ is chosen as in \eqref{eq_M} and satisfies $M>1$, and the shifts satisfy $\|[B_j^HB_j]^{-1}\|_2\leq1/S$, then \eqref{bound_rand} also holds.
\end{theorem}

\begin{proof}
The mean square approximation error for a given function $f \in \mathcal{K}_{d,\alpha, \bsgamma}$ is given by
\begin{align*}
\mathbb{E}_{\bsDelta} \left[ \|f - \mathcal{A}_{\bsDelta}(f) \|_2^2 \right] & = \sum_{\bsk \in A_{\alpha,\bsgamma,M}} \mathbb{E}_{\bsDelta} \left[ |\widehat{f}(\bsk) - \mathcal{B}_N(f, \bsk, \{ \bsy_m^{(s)}\}, \bsDelta )|^2 \right] + \sum_{\bsk \in \mathbb{Z}^d \setminus A_{\alpha,\bsgamma, M}} |\widehat{f}(\bsk)|^2.
\end{align*}

For the second term, we have
\begin{equation*}
\sum_{\bsk \in \mathbb{Z}^d \setminus A_{\alpha, \bsgamma,M}} |\widehat{f}(\bsk)|^2 = \sum_{\bsk \in \mathbb{Z}^d \setminus A_{\alpha, \bsgamma, M}} |\widehat{f}(\bsk)|^2 \frac{r^2_{\alpha,\bsgamma}(\bsk)}{r^2_{\alpha,\bsgamma}(\bsk)} \le \frac{1}{M^2} \sum_{\bsk \in \mathbb{Z}^d \setminus A_{\alpha, \bsgamma, M}} |\widehat{f}(\bsk)|^2 r^2_{\alpha,\bsgamma}(\bsk) \le \frac{\|f\|^2_{d,\alpha, \bsgamma} }{M^2}.
\end{equation*}

To analyse the first term, let $\bsell^{(1)}, \bsell^{(2)}, \ldots, \bsell^{(J)}$ be such that $\Gamma_{\alpha, \bsgamma,N}(\bsell^{(i)}; \bsg) \cap \Gamma_{\alpha, \bsgamma,N}(\bsell^{(j)}; \bsg) = \emptyset$ for all $i \neq j$ and $\bigcup_{j=1}^J \Gamma_{\alpha, \bsgamma, N}(\bsell^{(j)}; \bsg) = A_{\alpha,\bsgamma, M}$. Let $v_j = |\Gamma_{\alpha, \bsgamma, N}(\bsell^{(j)}; \bsg)|$. Then
\begin{align*}
\lefteqn{ \sum_{\bsk \in A_{\alpha, \bsgamma, M}} |\widehat{f}(\bsk) - \mathcal{B}_N(f, \bsk, \{\bsy_m^{(s)}\},   \bsDelta )|^2 } \\ & = \sum_{j=1}^J \sum_{\bsk \in \Gamma_{\alpha,\bsgamma, N}(\bsell^{(j)}; \bsg)} |\widehat{f}(\bsk) - \mathcal{B}_N(f, \bsk, \{\bsy_m^{(s)} \},  \bsDelta )|^2  \\ & = \sum_{j=1}^J \left\|  \begin{pmatrix}  \widehat{f}(\bsell_1) \\ \vdots \\ \widehat{f}(\bsell_{v_j}) \end{pmatrix}  - D(\bsDelta)^{-1} \left[B_j^H B_j \right]^{-1} B_j^{H}  \begin{pmatrix} \mathrm{e}^{2\pi \mathrm{i} \bsell_1 \cdot (\bsy_1^{(1)} + \bsDelta ) } \mathcal{F}_{N}(f,\bsell_1, \bsy_1^{(1)} + \bsDelta ) \\ \vdots \\ \mathrm{e}^{2\pi \mathrm{i} \bsell_{v_j} \cdot (\bsy_{v_j}^{(S)} + \bsDelta) } \mathcal{F}_{N}(f,\bsell_{v_j}, \bsy_{v_j}^{(S)} + \bsDelta) \end{pmatrix} \right\|_2^2 \\ & = \sum_{j=1}^J \left\| D(\bsDelta) \begin{pmatrix}  \widehat{f}(\bsell_1) \\ \vdots \\ \widehat{f}(\bsell_{v_j}) \end{pmatrix}  -  \left[B_j^H B_j \right]^{-1} B_j^{H}  \begin{pmatrix} \mathrm{e}^{2\pi \mathrm{i} \bsell_1 \cdot (\bsy_1^{(1)} + \bsDelta ) } \mathcal{F}_{N}(f,\bsell_1, \bsy_1^{(1)} + \bsDelta ) \\ \vdots \\ \mathrm{e}^{2\pi \mathrm{i} \bsell_{v_j} \cdot (\bsy_{v_j}^{(S)} + \bsDelta) } \mathcal{F}_{N}(f,\bsell_{v_j}, \bsy_{v_j}^{(S)} + \bsDelta) \end{pmatrix} \right\|_2^2.
\end{align*}
We have
\begin{align*}
\mathcal{F}_{N}(f, \bsell_m, \bsy_m^{(s)} + \bsDelta )  & =  \frac{1}{N} \sum_{n=0}^{N-1} f(\{n \bsg / N + \bsy_m^{(s)} + \bsDelta \}) \mathrm{e}^{-2\pi \mathrm{i} \bsell_m \cdot  (n \bsg / N + \bsy_m^{(s)} + \bsDelta ) }  \\ & = \sum_{\bsk - \bsell_m \in \mathcal{L}_N^\perp(\bsg)} \widehat{f}(\bsk) \mathrm{e}^{2\pi \mathrm{i} (\bsk - \bsell_m) \cdot ( \bsy_m^{(s)} + \bsDelta )}.
\end{align*}

Let $\{ \bsell_1, \ldots, \bsell_{v_j}\} = \Gamma_{\alpha,\bsgamma,N}(\bsell^{(j)};\bsg)$ where $v_j = R_{\alpha,\bsgamma,N}(\bsell^{(j)};\bsg)$. Then $\bsk-\bsell_m \in \mathcal{L}_N^\perp(\bsg)$ implies that $\bsk- \bsell^{(j)} \in \mathcal{L}_N^\perp(\bsg)$.  Then
\begin{align*}
& \left[ B_j^H B_j \right]^{-1}  B_j^H \left[ \begin{pmatrix} \mathrm{e}^{2\pi \mathrm{i} \bsell_1 \cdot (\bsy_1^{(1)} + \bsDelta) }\mathcal{F}_{N}(f,\bsell_1, \bsy_1^{(1)} + \bsDelta ) \\ \vdots \\ \mathrm{e}^{2\pi \mathrm{i} \bsell_{v_j} \cdot (\bsy_{v_j}^{(S)} + \bsDelta) } \mathcal{F}_{N}(f,\bsell_{v_j}, \bsy_{v_j}^{(S)} + \bsDelta ) \end{pmatrix} - B_j D(\bsDelta) \begin{pmatrix}  \widehat{f}(\bsell_1) \\ \vdots \\ \widehat{f}(\bsell_{v_j}) \end{pmatrix} \right]  \\  & = \left[ B_j^H B_j \right]^{-1}  B_j^H  \left( \sum_{\substack{ \bsk - \bsell_m \in \mathcal{L}_N^\perp(\bsg)  \\ \bsk \notin \Gamma_{\alpha, \bsgamma, N}(\bsell_m; \bsg) } } \widehat{f}(\bsk) \mathrm{e}^{2\pi \mathrm{i} \bsk   \cdot ( \bsy_m^{(s)} + \bsDelta) } \right)_{m = 1, \ldots, v_j, s = 1, \ldots, S} \\ & = \left[ B_j^H B_j \right]^{-1}  B_j^H \left( \sum_{\substack{ \bsk - \bsell^{(j)} \in \mathcal{L}_N^\perp(\bsg)  \\ \bsk \notin \Gamma_{\alpha, \bsgamma, N}(\bsell^{(j)}; \bsg) } } \widehat{f}(\bsk) \mathrm{e}^{2\pi \mathrm{i} \bsk   \cdot ( \bsy_m^{(s)} + \bsDelta ) } \right)_{m = 1, \ldots, v_j, s= 1, \ldots, S} \\ & =  \sum_{\substack{ \bsk - \bsell^{(j)} \in \mathcal{L}_N^\perp(\bsg)  \\ \bsk \notin \Gamma_{\alpha,\bsgamma, N}(\bsell^{(j)}; \bsg) } } \widehat{f}(\bsk) \left[ B_j^H B_j \right]^{-1}  B_j^H \underbrace{\begin{pmatrix} \mathrm{e}^{2\pi \mathrm{i} \bsk  \cdot ( \bsy_1^{(1)} + \bsDelta ) } \\ \vdots \\ \mathrm{e}^{2\pi \mathrm{i} \bsk \cdot ( \bsy_{v_j}^{(S)} + \bsDelta ) } \end{pmatrix} }_{=: \bsa_{\bsk}(\bsDelta)}.
\end{align*}
Thus
\begin{align*}
& \left\|  \begin{pmatrix}  \widehat{f}(\bsell_1) \\ \vdots \\ \widehat{f}(\bsell_{v_j}) \end{pmatrix} - D(\bsDelta)^{-1} \left[ B_j^H B_j \right]^{-1}  B_j^H \begin{pmatrix} \mathrm{e}^{2\pi \mathrm{i} \bsell_1 \cdot (\bsy_1^{(1)} + \bsDelta) } \mathcal{F}_{N}(f,\bsell_1, \bsy_1^{(1)} + \bsDelta ) \\ \vdots \\ \mathrm{e}^{2\pi \mathrm{i} \bsell_{v_j} \cdot (\bsy_{v_j}^{(S)} + \bsDelta) } \mathcal{F}_{N}(f,\bsell_{v_j}, \bsy_{v_j}^{(S)} + \bsDelta ) \end{pmatrix} \right\|_2^2 \\ & = \sum_{\substack{ \bsk - \bsell^{(j)}, \bsk' - \bsell^{(j)} \in \mathcal{L}_N^\perp(\bsg)  \\ \bsk, \bsk' \notin \Gamma_{\alpha, \bsgamma, N}(\bsell^{(j)}; \bsg) } } \widehat{f}(\bsk) \overline{\widehat{f}(\bsk')} c(\bsk, \bsk', \bsDelta), 
\end{align*}
where $c(\bsk, \bsk', \bsDelta) = e(\bsk, \bsDelta)^H e(\bsk', \bsDelta)$ and
\begin{align*}
e(\bsk, \bsDelta) & = \left[ B_j^H B_j \right]^{-1}  B_j^H \bsa_{\bsk}(\bsDelta)  = \mathrm{e}^{2\pi \mathrm{i} \bsk \cdot \bsDelta } \left[ B_j^H B_j \right]^{-1}  B_j^H \bsa_{\bsk}(\bszero).
\end{align*}
Hence
\begin{align*}
& \left\|  \begin{pmatrix}  \widehat{f}(\bsell_1) \\ \vdots \\ \widehat{f}(\bsell_{v_j}) \end{pmatrix} - D(\bsDelta)^{-1} \left[ B_j^H B_j \right]^{-1}  B_j^H \begin{pmatrix} \mathrm{e}^{2\pi \mathrm{i} \bsell_1 \cdot (\bsy_1^{(1)} + \bsDelta) } \mathcal{F}_{N}(f,\bsell_1, \bsy_1^{(1)} + \bsDelta ) \\ \vdots \\ \mathrm{e}^{2\pi \mathrm{i} \bsell_{v_j} \cdot (\bsy_{v_j}^{(S)} + \bsDelta) } \mathcal{F}_{N}(f,\bsell_{v_j}, \bsy_{v_j}^{(S)} + \bsDelta ) \end{pmatrix} \right\|_2^2 \\ & = \sum_{\substack{ \bsk - \bsell^{(j)}, \bsk' - \bsell^{(j)} \in \mathcal{L}_N^\perp(\bsg)  \\ \bsk, \bsk' \notin \Gamma_{\alpha, \bsgamma, N}(\bsell^{(j)}; \bsg) } } \widehat{f}(\bsk) \overline{\widehat{f}(\bsk')} \mathrm{e}^{2\pi \mathrm{i} (\bsk'-\bsk)\cdot \bsDelta}  \bsa_{\bsk}(\bszero)^H B_j [B_j^H B_j]^{-1} [B_j^H B_j]^{-1} B_j^H \bsa_{\bsk'}(\bszero).
\end{align*}

Since
\begin{equation*}
\mathbb{E}_{\bsDelta}[c(\bsk, \bsk', \bsDelta)] = \mathbb{E}_{\bsDelta}[\mathrm{e}^{2\pi \mathrm{i} (\bsk'-\bsk) \cdot \bsDelta} ] \bsa_{\bsk}(\bszero)^H B_j [B_j^H B_j]^{-1} [B_j^H B_j]^{-1} B_j^H \bsa_{\bsk'}(\bszero),
\end{equation*}
for $\bsk \neq \bsk' $ we get $\mathbb{E}_{\bsDelta}[c(\bsk, \bsk', \bsDelta)] = 0$ and
\begin{equation*}
\mathbb{E}_{\bsDelta}(c(\bsk, \bsk, \bsDelta)) =  \bsa_{\bsk}(\bszero)^H B_j [B_j^H B_j]^{-1} [B_j^H B_j]^{-1} B_j^H \bsa_{\bsk}(\bszero) =: c(\bsk).
\end{equation*}
Each component of $B_j^H\bsa_{\bsk}(\bszero)$ has modulus at most $v_jS$. Hence, using $\|[B_j^HB_j]^{-1}\|_2\leq1/S$,
\[
|c(\bsk)|
=\|(B_j^HB_j)^{-1}B_j^H\bsa_{\bsk}(\bszero)\|_2^2
\leq v_j^3
\leq R_{\alpha,\bsgamma,N}(\bsg)^3.
\]
Hence
\begin{align*}
& \mathbb{E}_{\bsDelta} \left\|  \begin{pmatrix}  \widehat{f}(\bsell_1) \\ \vdots \\ \widehat{f}(\bsell_{v_j}) \end{pmatrix} - D(\bsDelta)^{-1} \left[ B_j^H B_j \right]^{-1}  B_j^H \begin{pmatrix} \mathrm{e}^{2\pi \mathrm{i} \bsell_1 \cdot (\bsy_1^{(1)} + \bsDelta) } \mathcal{F}_{N}(f,\bsell_1, \bsy_1^{(1)} + \bsDelta ) \\ \vdots \\ \mathrm{e}^{2\pi \mathrm{i} \bsell_{v_j} \cdot (\bsy_{v_j}^{(S)} + \bsDelta) } \mathcal{F}_{N}(f,\bsell_{v_j}, \bsy_{v_j}^{(S)} + \bsDelta ) \end{pmatrix} \right\|_2^2 \\ & \le R_{\alpha,\bsgamma,N}(\bsg)^3 \sum_{\substack{\bsk - \bsell^{(j)} \in \mathcal{L}_N^\perp(\bsg)  \\ \bsk \notin \Gamma_{\alpha, \bsgamma, N}(\bsell^{(j)}; \bsg) } } |\widehat{f}(\bsk) |^2. 
\end{align*}
Thus
\begin{align*}
\sum_{\bsk \in A_{\alpha, \bsgamma, M}} \mathbb{E}_{\bsDelta} \left[ |\widehat{f}(\bsk) - \mathcal{B}_N(f, \bsk, \{\bsy_m^{(s)}\}, \bsDelta)|^2 \right] & \le R_{\alpha,\bsgamma,N}(\bsg)^3 \sum_{j=1}^J \sum_{\substack{\bsk - \bsell^{(j)} \in \mathcal{L}_N^\perp(\bsg)  \\ \bsk \notin \Gamma_{\alpha, \bsgamma, N}(\bsell^{(j)}; \bsg) } } |\widehat{f}(\bsk) |^2 \\ & \le R_{\alpha,\bsgamma,N}(\bsg)^3 \sum_{\bsk \in \mathbb{Z}^d \setminus A_{\alpha, \bsgamma, M}} |\widehat{f}(\bsk)|^2 \frac{r_{\alpha,\bsgamma}^2(\bsk)}{r_{\alpha,\bsgamma}^2(\bsk)} \\ & \le R_{\alpha,\bsgamma,N}(\bsg)^3 \|f\|_{d, \alpha, \bsgamma}^2 \max_{\bsk \in \mathbb{Z}^d \setminus A_{\alpha, \bsgamma, M}} \frac{1}{r_{\alpha,\bsgamma}^2(\bsk)} \\ & \le R_{\alpha,\bsgamma,N}(\bsg)^3 \frac{\|f\|_{d, \alpha, \bsgamma}^2}{M^2}.
\end{align*}

Thus we obtain
\begin{equation*}
\mathbb{E}_{\bsDelta} \left[ \|f - \mathcal{A}_{\bsDelta}(f) \|_2^2 \right] \le \frac{\|f\|_{d, \alpha, \bsgamma}^2}{M^2} \left(1 + R_{\alpha,\bsgamma,N}(\bsg)^3\right).
\end{equation*}
Since the bound holds for all functions $f$ in $\mathcal{K}_{d,\alpha, \bsgamma}$, we have
\begin{equation*}
\mathrm{e}_2^{\ran}((\mathcal{A}_{\bsDelta})_{\bsDelta},\mathcal{K}_{d,\alpha, \bsgamma}) \le \frac{\sqrt{ 1+ R_{\alpha,\bsgamma,N}(\bsg)^3 } }{M}.
\end{equation*}
Fix $\lambda\in(1/\alpha,2]$. From \eqref{eq:MN2}, $M\geq c_{\lambda,d,\alpha,\bm\gamma,\delta}N^{1/\lambda}$, and \eqref{bound_g_exist} therefore yields
\[
\mathrm e_2^{\ran}((\mathcal A_{\bm\Delta})_{\bm\Delta},\mathcal K_{d,\alpha,\bm\gamma})
\leq C_{\lambda,d,\alpha,\bm\gamma,\delta}(\log N)^{3(d-1)/2}N^{-1/\lambda}.
\]
Using \eqref{bound_Np}, this is bounded by $C p^{-1/\lambda}(\log p)^b$ for a constant $b>0$. Given $\tau>0$, choose $\lambda$ sufficiently close to $1/\alpha$ that $1/\lambda>\alpha-\tau/2$, and absorb $(\log p)^b$ into $p^{\tau/2}$. This proves \eqref{bound_rand}. 

The result for the case where $\bm{g}$ is constructed by Algorithm~\ref{cbc_principle} follows by similar arguments.
\end{proof}

\section{Numerical experiments}\label{sec:numerical_experiments}

\subsection{Aliasing}\label{subsec:aliasing}

We first test the maximum number of aliased frequencies $R_{\alpha, \bm{\gamma}, N}(\bsg)$ for random choices of the generating vector $\bm{g}$. 

\begin{figure}[htbp]
  \centering
  \begin{minipage}[b]{0.48\textwidth}
    \centering
    \includegraphics[width=\textwidth]{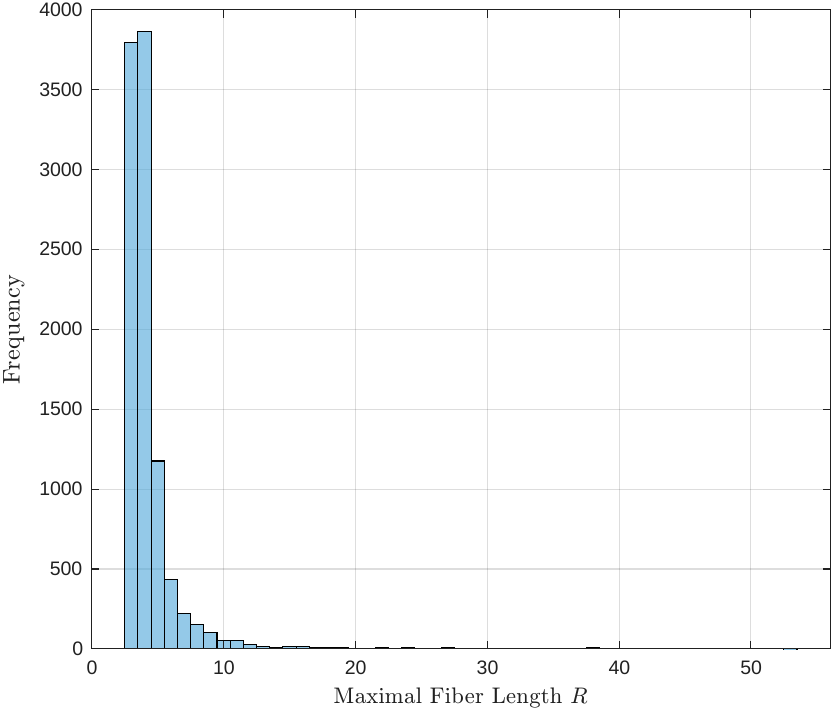}
    \subcaption{$d=2$, $N=9973$, $M=355$}
    \label{fig:d2_Nsmall}
  \end{minipage}
  \hfill
  \begin{minipage}[b]{0.48\textwidth}
    \centering
    \includegraphics[width=\textwidth]{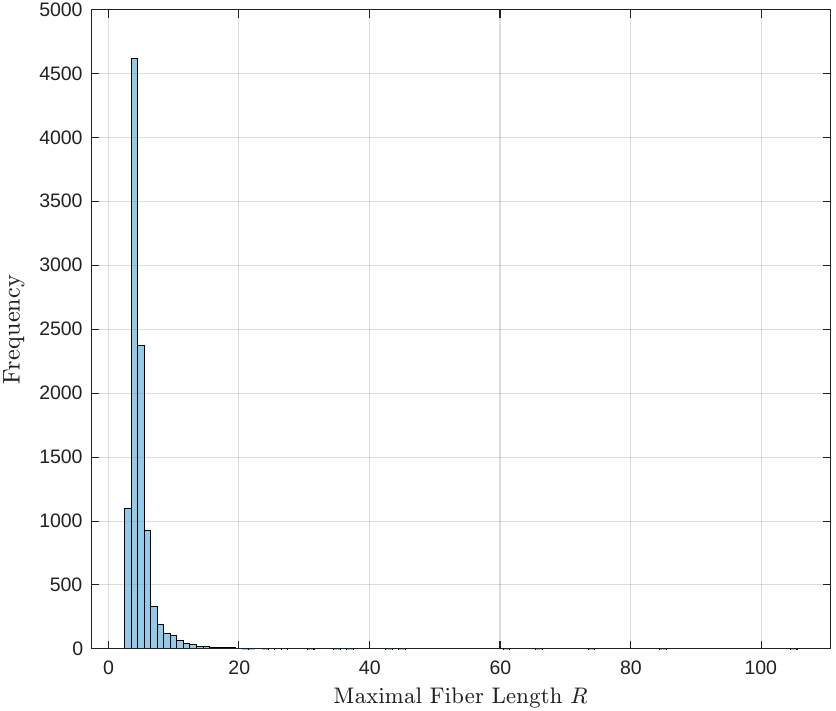}
    \subcaption{$d=2$, $N=99991$, $M=2754$}
    \label{fig:d2_Nlarge}
  \end{minipage}

  \vspace{1em}

  \begin{minipage}[b]{0.48\textwidth}
    \centering
    \includegraphics[width=\textwidth]{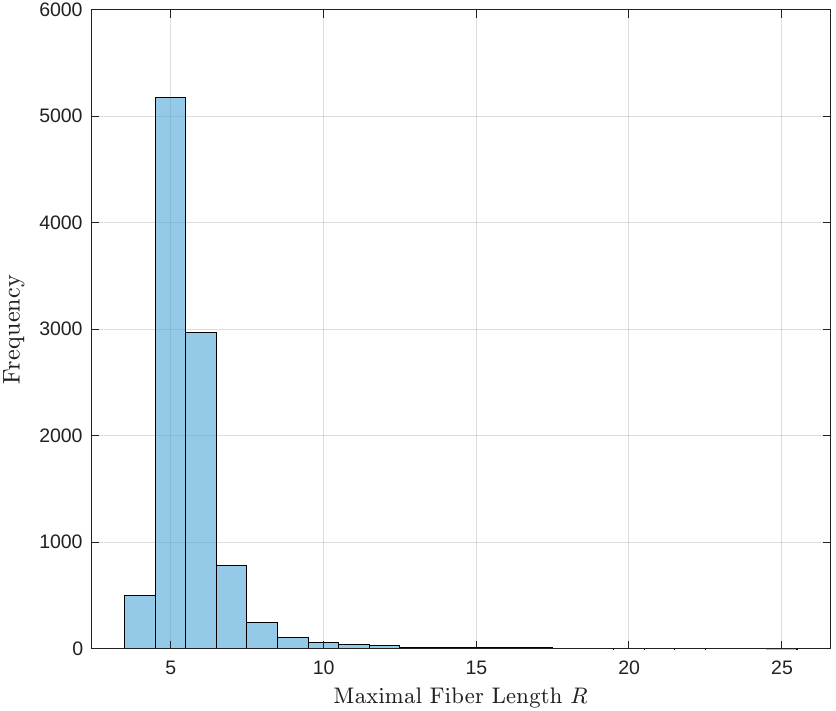}
    \subcaption{$d=100$, $N=9973$, $M=598$}
    \label{fig:d100_Nsmall}
  \end{minipage}
  \hfill
  \begin{minipage}[b]{0.48\textwidth}
    \centering
    \includegraphics[width=\textwidth]{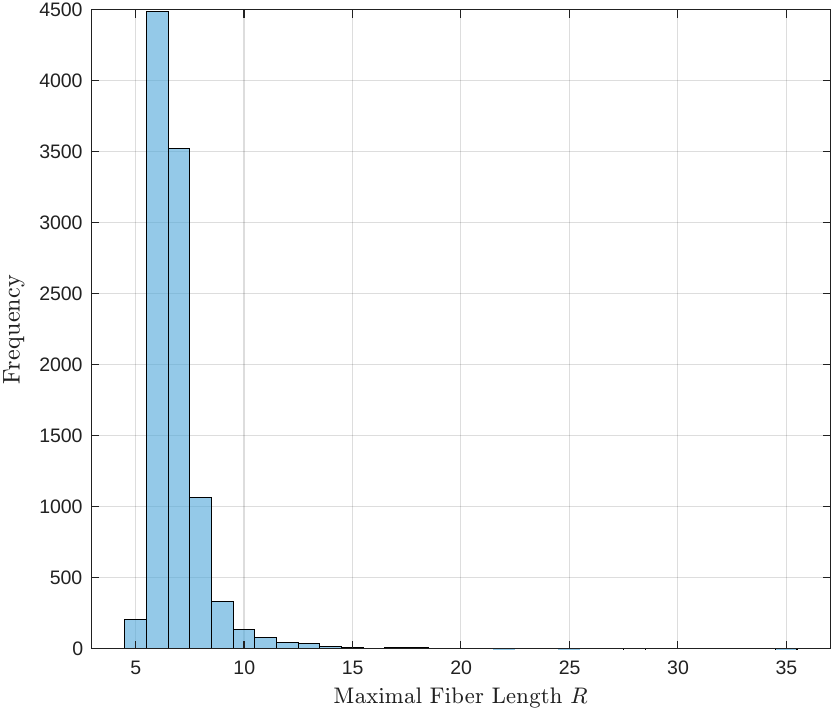}
    \subcaption{$d=100$, $N=99991$, $M=3600$}
    \label{fig:d100_Nlarge}
  \end{minipage}

  \caption{Histograms of the maximal fiber length $R$ for $10^4$ random choices of the generating vector $\boldsymbol{g}$. The top row (a, b) shows the unweighted case for $d=2$, and the bottom row (c, d) shows the weighted case for $d=100$ with weights $\gamma_j = j^{-2}$. The smoothness parameter $\alpha$ is set to $1$ in all cases. In each case, the parameter $M$ is chosen such that the cardinality of the frequency index set satisfies $|A_{\alpha,\boldsymbol{\gamma},M}| \approx N$. In particular, $M=355, 2754, 598$, and $3600$, respectively.}
  \label{fig:hist2}
\end{figure}

In Figure~\ref{fig:hist2}, for a fixed smoothness $\alpha=1$, we empirically evaluate the distribution of the maximal fiber length $R=R_{\alpha, \bm{\gamma}, N}(\bsg)$ across four distinct configurations, comparing the unweighted, low-dimensional case ($d=2$, $\gamma_1=\gamma_2=1$) and the weighted, high-dimensional case ($d=100$, $\gamma_j=j^{-2}$) for two different prime numbers $N \approx 10^4$ and $N \approx 10^5$. 
In each setting, the parameter $M$ is chosen so that the cardinality of the frequency index set satisfies $|A_{\alpha,\boldsymbol{\gamma},M}| \approx N$, and the generating vectors are drawn independently and randomly from the set $\{1,\ldots,N-1\}^d$ for a total of $10^4$ times. Specifically, for $N = 9973$, we set $M = 355$ ($|A| = 9969$) for $d=2$, and $M = 598$ ($|A| = 9973$) for $d=100$. For $N = 99991$, we choose $M = 2754$ ($|A| = 99925$) for $d=2$, and $M = 3600$ ($|A| = 98983$) for $d=100$. The case $|A|=N$ lies just outside the strict cardinality assumption used in Lemma~\ref{lem_g_exist}; it is included here as a numerical stress test.

The statistical properties of $R$ over $10^4$ randomly chosen generating vectors are summarized as follows:
\begin{itemize}
    \item For $(d,N)=(2,9973)$: its minimum, median, and maximum are $3$, $4$, and $53$.
    \item For $(d,N)=(2,99991)$: its minimum, median, and maximum are $3$, $4$, and $105$.
    \item For $(d,N)=(100,9973)$: its minimum, median, and maximum are $4$, $5$, and $25$.
    \item For $(d,N)=(100,99991)$: its minimum, median, and maximum are $5$, $7$, and $35$.
\end{itemize}

The main purpose of this experiment is to illustrate the size of the fibers that occur in practice. While the theoretical bounds are formulated in terms of the worst-case quantity $R_{\alpha,\bm{\gamma},N}(\bm{g})$, randomly chosen generating vectors often produce fibers whose maximal size is much smaller than what these bounds might suggest. The available worst-case upper bound in Lemma~\ref{lem_g_exist} has exponential dependence on the dimension $d$, irrespective of how rapidly the weights decay. The empirical result, observed above, has a direct impact on the efficiency of the method in the weighted setting, since the least-squares reconstruction of the aliased Fourier coefficients is performed separately on each fiber. If a fiber has cardinality $v$, then the corresponding linear system involves $v$ unknown Fourier coefficients, and the computational effort of this de-aliasing step grows at least quadratically with $v$. The small values of $R_{\alpha,\bm{\gamma},N}(\bm{g})$ observed in the numerical experiment suggest that the actual computational effort is significantly smaller than what the theoretical bounds predict.

\subsection{Approximation}

\subsubsection*{Low-dimensional examples}

In the numerical approximation experiments below, the radius $M$ is selected heuristically to assess practical performance and is not intended to satisfy every sufficient parameter condition in Theorems~\ref{thm_det} and~\ref{thm:L_2_rce}; in particular, the rows with $d=3$ and $N=131$ have $|A_{\alpha,\bsgamma,M}|=135>N$. We use $S=\lceil2KR\log N\rceil$, as in Lemma~\ref{lem_eig}, and accept a sampled set of shifts only if
\[
\max_{1\leq j\leq J}S\left\|(B_j^HB_j)^{-1}\right\|_2\leq1;
\]
otherwise, the shifts are resampled. Thus the numerical experiments use the required conditioning property directly rather than relying on the conservative union-bound probability estimate.

First, let us consider the deterministic algorithm, i.e., we do not apply a random shift to the algorithm. Let $d=2$ and the target function be given by $f(x,y) = \exp(\cos(2\pi x)+\sin(2 \pi y))$. Notice that this function has fast-decaying Fourier coefficients. We choose a set of prime numbers $N$, randomly select a generating vector $\bm{g} = (1,g)$ with $g \in \{1, 2, \ldots, N-1\}$ and random vectors $\{\bm{y}_m^{(s)}\}$. We fix $K=1.1$, and $M$ is heuristically set to $N^{0.85}/7.35$ so that $|A_{\alpha,\bsgamma,M}|<N$. The result is shown in Table~\ref{table1}. To assess the numerical performance of the deterministic algorithm, we evaluate the $L_{\infty}$ error of the algorithm based on a large set of random points in $[0,1]^2$.

\begin{table}[tbp]
\begin{tabular}{rrrrrr}
\hline
$N$ & $g$ & $S$ & $|A_{\alpha, \bm{\gamma}, M}|$ & $R_{\alpha, \bm{\gamma}, N}$ & error \\
\hline
19   & 11  & 7  & 9    & 1 & $1.584\times 10^{0}$  \\
53   & 6  & 18  & 33   & 2 & $1.303\times 10^{-1}$ \\
131  & 127 & 43 & 113  & 4 & $2.851\times 10^{-3}$ \\
311  & 292 & 38 & 277  & 3 & $1.010\times 10^{-5}$ \\
719  & 498 & 44 & 705  & 3 & $5.692\times 10^{-10}$ \\
1619 & 1163 & 49 & 1593 & 3 & $2.095\times 10^{-14}$ \\
\hline
\end{tabular}
\caption{Numerical experiment using the deterministic algorithm. This table shows the number of points of the lattice rule $N$ with generating vector $(1,g)$, where $g$ is shown in the second column. The remaining columns show the chosen $S$, the size of the set of frequencies $A_{\alpha, \bm{\gamma}, M}$, the maximal fiber length $R_{\alpha, \bm{\gamma}, N}$, and the approximation error in $L_\infty$ norm (which we approximated using random points). In this experiment $\alpha = 1$ and $\bm{\gamma} = (1,1)$. }\label{table1}
\end{table}


As another example, we consider the case $d=3$ and the function $f(x,y,z) = \exp(\cos(2\pi x)+\sin(2 \pi y)) B_2(z)$, where $B_2(z) = z^2-z+1/6$ is the second Bernoulli polynomial. The Fourier coefficients for the second Bernoulli polynomial are known to be of order $|k|^{-2}$. 
We choose a set of prime numbers $N$, randomly select a generating vector $\bm{g} = (1,g_2, g_3)$ with $g_2, g_3 \in \{1, 2, \ldots, N-1\}$ and random vectors $\{\bm{y}_m^{(s)}\}$. We fix $K=1.1$ and $M$ is heuristically set to $N^{0.65}/7.35$. The result is shown in Table~\ref{table2}. Here again, the $L_{\infty}$ error is evaluated.

\begin{table}[tbp]
\begin{tabular}{rrrrrrr}
\hline
$N$ & $g_2$ & $g_3$ & $S$ & $|A_{\alpha, \bm{\gamma}, M}|$ & $R_{\alpha, \bm{\gamma},N}$ & error \\
\hline
53   & 6   & 45   & 18 & 27   & 2 & $5.597\times 10^{-1}$ \\
131  & 47   & 82   & 33 & 135  & 3 & $1.984\times 10^{-1}$ \\
311  & 187   & 59  & 64 & 279  & 5 & $1.133\times 10^{-1}$ \\
719  & 630  & 339  & 44 & 683  & 3 & $8.797\times 10^{-2}$ \\
1619 & 722  & 1394  & 66 & 1577 & 4 & $5.425\times 10^{-2}$ \\
3671 & 3445 & 483 & 91 & 3349 & 5 & $2.687\times 10^{-2}$ \\
8161 & 1267 & 6939 & 100 & 6499 & 5 & $1.610\times 10^{-2}$ \\
\hline
\end{tabular}
\caption{Numerical experiment using the deterministic algorithm. This table shows the number of points of the lattice rule $N$ with generating vector $(1,g_2, g_3)$, where $g_2$ is shown in the second column, $g_3$ in the third column. The remaining columns show the chosen $S$, the size of the set of frequencies $A_{\alpha, \bm{\gamma}, M}$, the maximal fiber length $R_{\alpha, \bm{\gamma}, N}$, and the approximation error in $L_\infty$ norm (which we approximated using random points). In this experiment $\alpha = 1$ and $\bm{\gamma} = (1,1,1)$. }\label{table2}
\end{table}


Next, let us consider the randomized algorithm, i.e., we apply a random shift to the algorithm. We test the same two functions as above. For the bivariate function, we choose a set of prime numbers $N$, randomly select a generating vector $\bm{g} = (1,g)$ with $g \in \{1, 2, \ldots, N-1\}$, and random vectors $\{\bm{y}_m^{(s)}\}$. For $10$ independent random shifts $\bm\Delta_1,\ldots,\bm\Delta_{10}$, we report the root-mean-square of the estimated $L_2$ errors, namely
\[
\left(\frac{1}{10}\sum_{q=1}^{10}\left\|f-\mathcal A_{\bm\Delta_q}(f)\right\|_2^2\right)^{1/2}.
\]
This is the sample analogue of the randomized error considered in Theorem~\ref{thm:L_2_rce}. The result is shown in Table~\ref{table4}. In this case, the corresponding theoretical result (Theorem~\ref{thm:L_2_rce}) concerns the randomized $L_2$ error, so we evaluate the $L_2$ error of the algorithm based on a large set of random points in $[0,1]^2$.

Likewise, for the trivariate function, we choose a set of prime numbers $N$, randomly select a generating vector $\bm{g} = (1,g_2, g_3)$ with $g_2, g_3 \in \{1, 2, \ldots, N-1\}$ and random vectors $\{\bm{y}_m^{(s)}\}$. The errors are aggregated using the same root-mean-square formula over $10$ independent random shifts. The result is shown in Table~\ref{table5}.

\begin{table}[tbp]
\begin{tabular}{rrrrrr}
\hline
$N$ & $g$ & $S$ & $|A_{\alpha, \bm{\gamma}, M}|$ & $R_{\alpha, \bm{\gamma}, N}$ & error \\
\hline
19   & 11  & 7  & 9    & 1 & $4.127\times 10^{-1}$  \\
53   & 3  & 27  & 33   & 3 & $3.736\times 10^{-2}$ \\
131  & 95  & 22 & 113  & 2 & $9.656\times 10^{-4}$ \\
311  & 166 & 26 & 277  & 2 & $2.513\times 10^{-6}$ \\
719  & 533 & 44 & 705  & 3 & $1.259\times 10^{-10}$ \\
1619 & 549 & 66 & 1593 & 4 & $5.233\times 10^{-15}$ \\
\hline
\end{tabular}
\caption{Numerical experiment using the randomized algorithm, reporting the root-mean-square over $10$ independent random shifts. This table shows the number of points of the lattice rule $N$ with generating vector $(1,g)$, where $g$ is shown in the second column. The remaining columns show the chosen $S$, the size of the set of frequencies $A_{\alpha, \bm{\gamma}, M}$, the maximal fiber length $R_{\alpha, \bm{\gamma}, N}$, and the approximation error in $L_2$ norm (which we approximated using random points). In this experiment $\alpha = 1$ and $\bm{\gamma} = (1,1)$.}\label{table4}
\end{table}




\begin{table}[tbp]
\begin{tabular}{rrrrrrr}
\hline
$N$ & $g_2$ & $g_3$ & $S$ & $|A_{\alpha,\bsgamma,M}|$ & $R_{\alpha,\bsgamma,N}$ & error \\
\hline
53   & 25   & 13   & 18 & 27   & 2 & $7.388\times 10^{-2}$ \\
131  & 92   & 89   & 33 & 135  & 3 & $1.660\times 10^{-2}$ \\
311  & 45  & 129  & 38 & 279  & 3 & $8.166\times 10^{-3}$ \\
719  & 107  & 421  & 44 & 683  & 3 & $3.103\times 10^{-3}$ \\
1619 & 1510  & 61  & 49 & 1577 & 3 & $1.597\times 10^{-3}$ \\
3671 & 1752 & 2645 & 73 & 3349 & 4 & $6.912\times 10^{-4}$ \\
8161 & 4900 & 7128 & 60 & 6499 & 3 & $3.119\times 10^{-4}$ \\
\hline
\end{tabular} 
\caption{Numerical experiment using the randomized algorithm, reporting the root-mean-square over $10$ independent random shifts. This table shows the number of points of the lattice rule $N$ with generating vector $(1,g_2, g_3)$, where $g_2$ is shown in the second column, $g_3$ in the third column. The remaining columns show the chosen $S$, the size of the set of frequencies $A_{\alpha, \bm{\gamma}, M}$, the maximal fiber length $R_{\alpha, \bm{\gamma}, N}$, and the approximation error in $L_2$ norm (which we approximated using random points). In this experiment $\alpha = 1$ and $\bm{\gamma} = (1,1,1)$. }\label{table5}
\end{table}

\subsubsection*{Comparison with existing algorithms}


We compare the empirical performance of the proposed multiple-shift lattice algorithm with several existing algorithms in a low-dimensional setting.

As a first test problem, we consider the bivariate kink function $$f(x,y) = \frac{5^{3/2} 15^2}{48} \max\left\{0, \frac{1}{5} - \left(x-\frac{1}{2}\right)^2\right\} \max\left\{0, \frac{1}{5} - \left(y-\frac{1}{2}\right)^2\right\},$$ which was also considered in \cite{BGKS25}. For our algorithm, we choose a set of prime numbers $N$. For each $N$, a generating vector $\bm{g} = (1,g)$ is selected uniformly at random from $g \in \{1, 2, \ldots, N-1\}$, and the random vectors $\{\bm{y}_m^{(s)}\}$ are generated without applying a random shift. To ensure consistency with the experimental setup in \cite{BGKS25}, we consider the weighted hyperbolic cross $A_{\alpha,\bsgamma,M}$ with $\alpha=1$ and $\gamma_1=\gamma_2=0.5$, where the threshold $M$ is chosen as large as possible such that $|A_{\alpha,\bsgamma,M}|<N$. 

For comparison with the subsampling method combined with least-squares approximation, we digitally extracted the error convergence from the upper panel of \cite[Figure~3]{BGKS25}. We also include results for the existing single rank-1 lattice algorithm satisfying the reconstruction property over the same hyperbolic cross $A_{\alpha,\bsgamma,M}$ \cite{Kammerer2013SJNA}. In this case, the number of points $N$ and the generating vector $\bsg$ are chosen such that the maximum fiber length over $A_{\alpha,\bsgamma,M}$ is equal to $1$. Furthermore, we tested the median lattice-based algorithm \cite[Algorithm~1]{PKG25} without a random shift but with the same parameter choice $\alpha=1$ and $\gamma_1=\gamma_2=0.5$. The remaining parameters were selected according to \cite[Remark~16]{PKG25}.

\begin{figure}[htbp]
  \centering
  \includegraphics[width=0.6\textwidth]{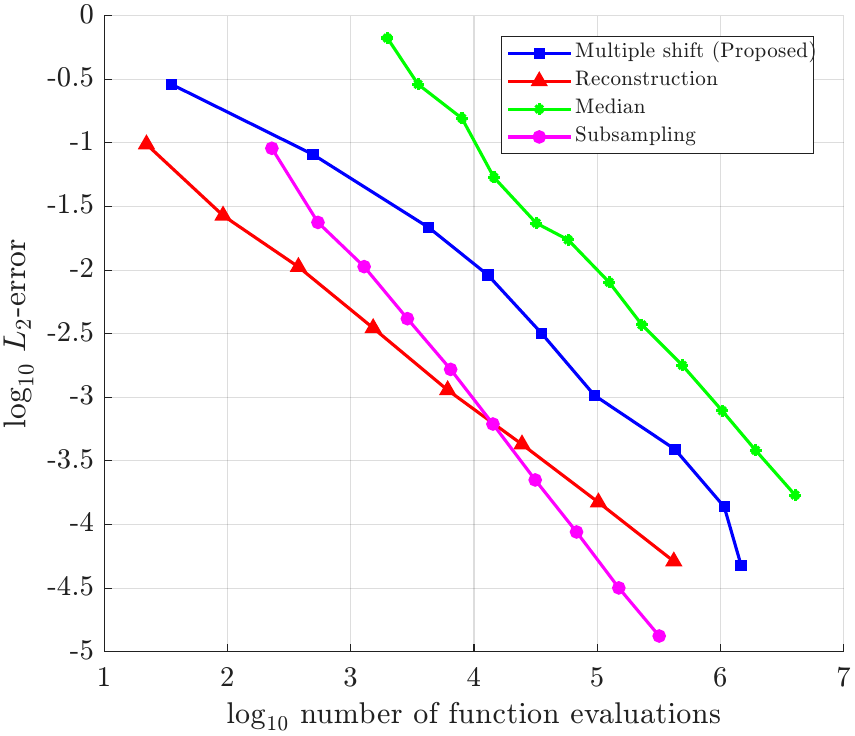}\\
\caption{Comparison of the $L_2$ approximation error of four different algorithms for the bivariate kink function. The $x$-axis shows the total number of function evaluations and the $y$-axis shows the $L_2$ approximation error.} \label{fig:comparison_d2}
\end{figure}

The results are shown in Figure~\ref{fig:comparison_d2}. The horizontal axis represents the total number of function evaluations, while the vertical axis represents the $L_2$ error. For each algorithm, the $L_2$ error is estimated by taking the square root of the average of $|f(\bsx)-\mathcal{A}(f)(\bsx)|^2$ over the first $2^{15}$ points of the two-dimensional Sobol' sequence. This deterministic evaluation procedure is adopted to ensure a fair comparison and to avoid unnecessary fluctuations caused by random sampling in the error estimation. As can be seen from Figure~\ref{fig:comparison_d2}, the reconstruction rank-1 lattice algorithm performs best when the number of function evaluations is relatively small. As the number of function evaluations increases, however, the subsampling algorithm eventually becomes the most accurate among the methods considered. Although the proposed multiple-shift algorithm does not outperform the reconstruction rank-1 lattice algorithm over the range of $N$ considered here, it exhibits a steeper empirical convergence rate, suggesting that it may become more competitive for larger values of $N$. The median lattice-based algorithm exhibits a similar trend.


\begin{figure}[htbp]
  \centering
  \includegraphics[width=0.6\textwidth]{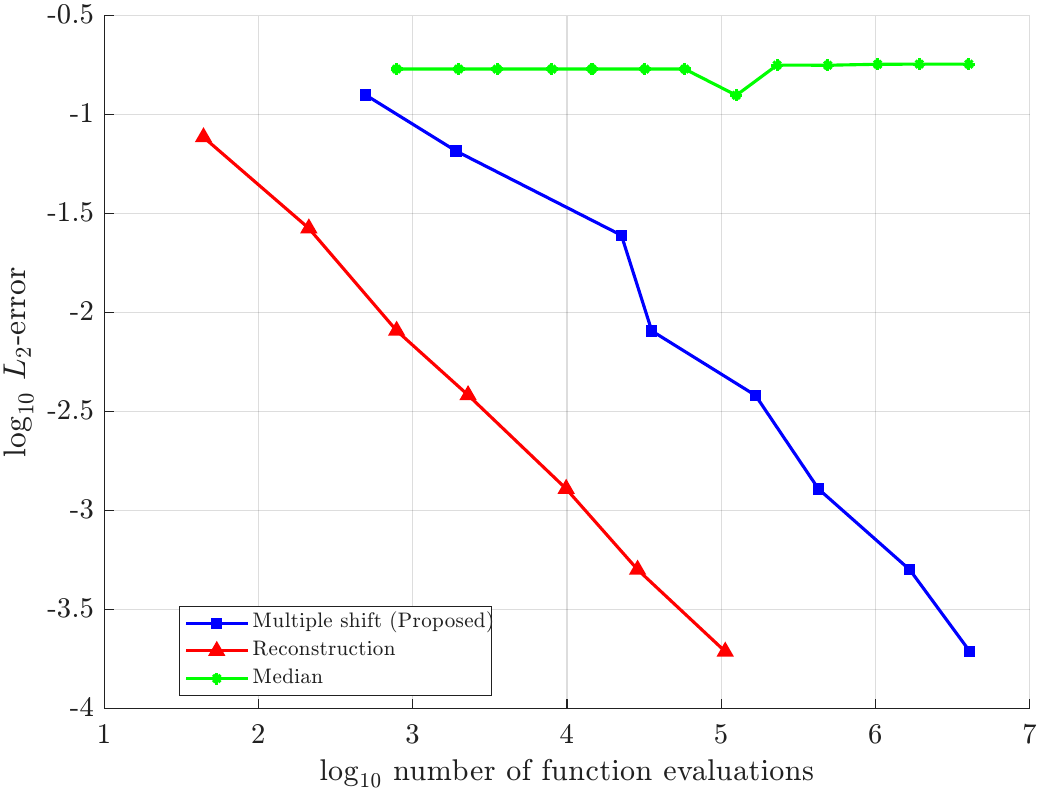}\\
\caption{Comparison of the $L_2$ approximation error of three different algorithms for the trivariate test function. The $x$-axis shows the total number of function evaluations and the $y$-axis shows the $L_2$ approximation error.} \label{fig:comparison_d3}
\end{figure}

We next consider the trivariate test function $f(x,y,z) = \exp(\cos(2\pi x)+\sin(2 \pi y)) B_2(z)$.  Here, we compare the proposed multiple-shift algorithm with the reconstruction rank-1 lattice algorithm and the median lattice-based algorithm. As in the previous example, we choose $\alpha=1$ and $\gamma_1=\gamma_2=\gamma_3=0.5$. The remaining parameters are selected in the same manner. The results are shown in Figure~\ref{fig:comparison_d3}. As in the two-dimensional example, the reconstruction rank-1 lattice algorithm performs best over the range of $N$ considered here. Although the proposed multiple-shift algorithm has a better asymptotic convergence rate in theory, this behavior is not yet evident from the numerical results, possibly because the range of $N$ is still too limited for the asymptotic regime to emerge. A similar observation applies to the median lattice-based algorithm, for which no clear error decay is observed over the range considered. This behavior is partly due to the parameter selection rule. That is, the choices of the threshold $M$ and the number of repetitions $R$ in \cite[Remark~16]{PKG25} are theoretically optimized for settings where the weights decay sufficiently fast. Thus, in our low-dimensional setting with nondecaying weights, the asymptotic convergence of the median-based algorithm is expected to kick in only for much larger values of $N$ than those tested here.

In the present comparison, we do not include multiple rank-1 lattice algorithms \cite{Kammerer2018JFAA,KV19}. A more comprehensive comparison among these distinct algorithms lies beyond the scope of this paper and is left for future research.

\subsubsection*{High-dimensional weighted examples}

Finally, we consider high-dimensional weighted examples. Although our theoretical results do not establish tractability of the algorithm, we empirically evaluate whether our algorithm effectively works in high dimensions. The test function we use is given by
\[ f(\bsx) = \prod_{j=1}^{d}\left[ 1+\omega_j \frac{121\sqrt{33}}{100} \max\left( \frac{25}{121}-\left( x_j-\frac{1}{2}\right)^2, 0 \right)\right], \]
where $\omega_j$ controls the relative importance of the variable $x_j$. In what follows, we consider the cases $d\in \{10, 20, 100\}$ and $\omega_j=j^{-2}$ or $j^{-4}$.

\begin{figure}[htbp]
  \centering
  \begin{minipage}[b]{0.48\textwidth}
    \centering
    \includegraphics[width=\textwidth]{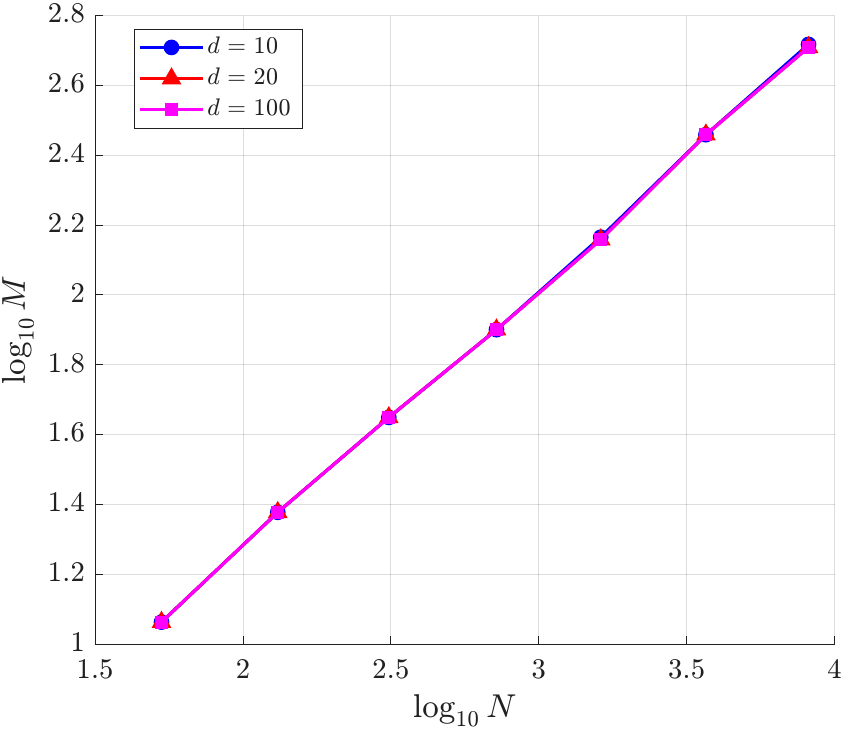}
    \subcaption{$\omega_j=j^{-2}$}
    \label{fig:size_slow}
  \end{minipage}
  \hfill
  \begin{minipage}[b]{0.48\textwidth}
    \centering
    \includegraphics[width=\textwidth]{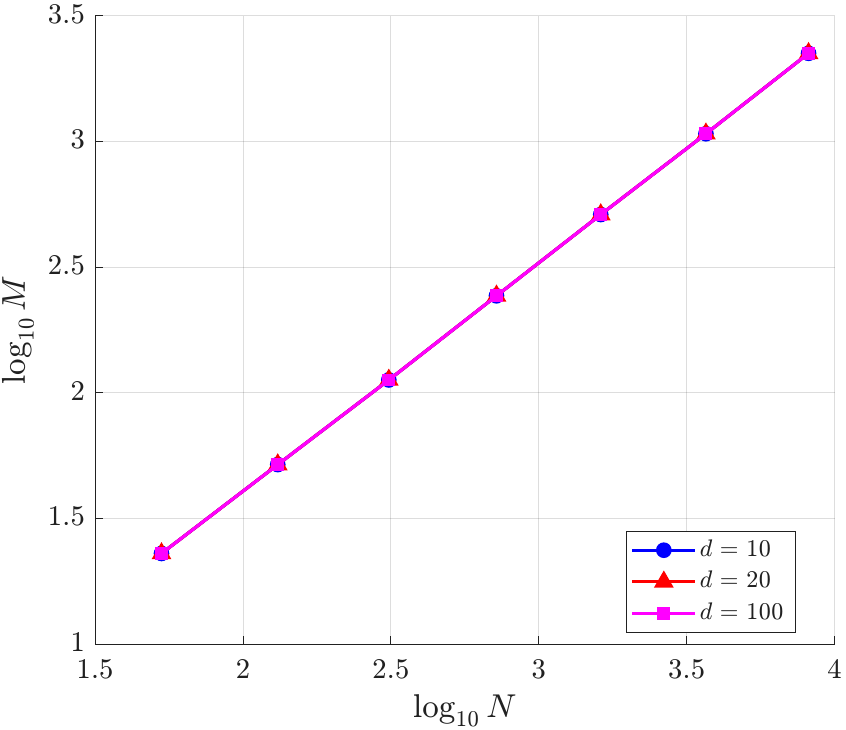}
    \subcaption{$\omega_j=j^{-4}$}
    \label{fig:size_fast}
  \end{minipage}
\caption{The radius of the weighted hyperbolic cross as a function of $N$ for the high-dimensional cases. The $x$-axis shows the number of points $N$ for a single rank-1 lattice rule, and the $y$-axis shows the radius $M$ of the weighted hyperbolic cross, which is made as large as possible so that $|A_{\alpha,\bsgamma,M}|<N$.} \label{fig:radius_hyperbolic}
\end{figure}

We choose a set of prime numbers $N\in \{ 53, 131, 311, 719, 1619, 3671, 8161\}$ (as done before), and randomly select a generating vector $\bm{g}$ with $g_1=1$ and $g_2,\ldots,g_d \in \{1, 2, \ldots, N-1\}$, and random vectors $\{\bm{y}_m^{(s)}\}$. We fix $K=1.1$. To construct the weighted hyperbolic cross, we employ $\alpha=1, \gamma_j=\omega_j$, and choose $M$ as large as possible such that $|A_{\alpha,\bsgamma,M}|<N$. In Figure~\ref{fig:radius_hyperbolic}, we show how $\log_{10} M$ grows as a function of  $\log_{10} N$. Using linear regression, we see that $M$ is of order $N^{0.75}$ for the case $\omega_j=j^{-2}$, and of order $N^{0.91}$ for the case $\omega_j=j^{-4}$. These observations align well with Lemma~\ref{lem:hyp_cross_size} and \eqref{eq:MN}. 

\begin{figure}[htbp]
  \centering
  \begin{minipage}[b]{0.48\textwidth}
    \centering
    \includegraphics[width=\textwidth]{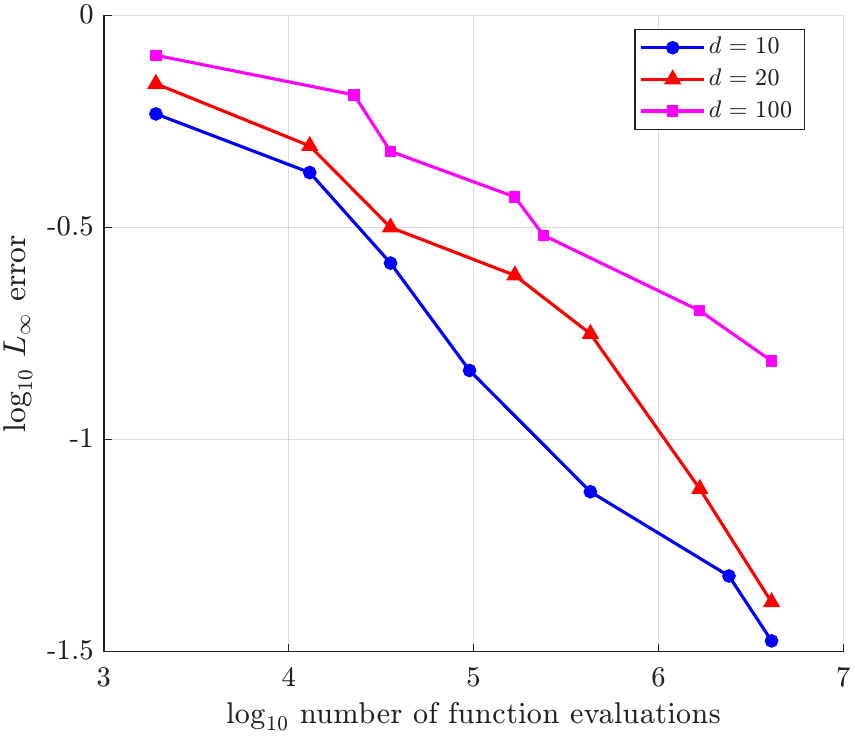}
    \subcaption{$\omega_j=j^{-2}$}
    \label{fig:high-dim_slow}
  \end{minipage}
  \hfill
  \begin{minipage}[b]{0.48\textwidth}
    \centering
    \includegraphics[width=\textwidth]{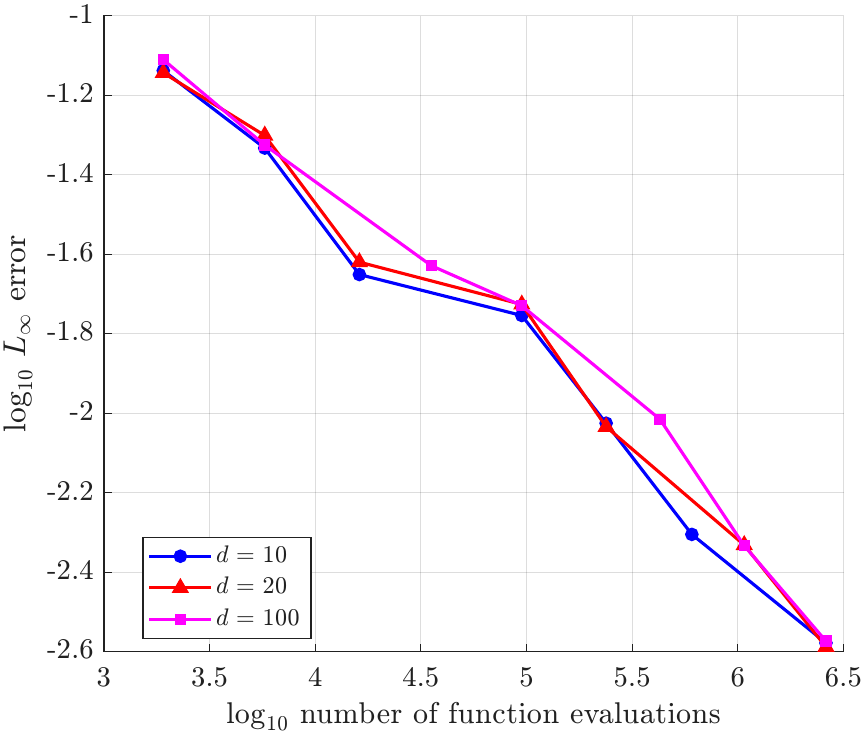}
    \subcaption{$\omega_j=j^{-4}$}
    \label{fig:high-dim_fast}
  \end{minipage}
\caption{Numerical experiment using the deterministic algorithm for high-dimensional weighted functions. The $x$-axis shows the total number of points used in the approximation algorithm and the $y$-axis shows the $L_{\infty}$ approximation error.} \label{fig:high-dim}
\end{figure}

The $L_\infty$ approximation error is shown in Figure~\ref{fig:high-dim}, where the left panel is for the case $\omega_j=j^{-2}$ and the right panel is for the case $\omega_j=j^{-4}$. The horizontal axis represents the total number of function evaluations, $p=NR_{\alpha,\bsgamma,N}(\bsg)S$, whereas the vertical axis represents the $L_{\infty}$ error. Here, we approximate the $L_{\infty}$ error of the algorithm by the maximum of the pointwise evaluations $|f(\bsx)-\mathcal{A}(f)(\bsx)|$ over the first $2^{15}$ points of the $d$-dimensional Sobol' sequence. Although we set the same numbers for $N$ for all dimensions $d$, the corresponding $p$ does not differ much. This is consistent with the empirical result about the maximum fiber length observed in Section~\ref{subsec:aliasing}; even in high dimensions, $R_{\alpha,\bsgamma,N}(\bsg)$ can remain small in the weighted setting. As a result, in particular, for the case $\omega_j=j^{-4}$, the weights decay fast enough that the error convergence does not differ much across dimensions. The observed empirical convergence rate with respect to $p$ is approximately $0.44$, while it reaches around $0.66$ when evaluated against $N$. Although these rates are inferior to the theoretical expectation implied by the smoothness of $f$, this discrepancy can be partially attributed to the sample size effect; for sufficiently large $N$ (and consequently larger $p$), the convergence rate is expected to improve. Crucially, however, the observed rates remain remarkably stable and do not deteriorate even in high dimensions.

\begin{figure}[htbp]
  \centering
  \begin{minipage}[b]{0.48\textwidth}
    \centering
    \includegraphics[width=\textwidth]{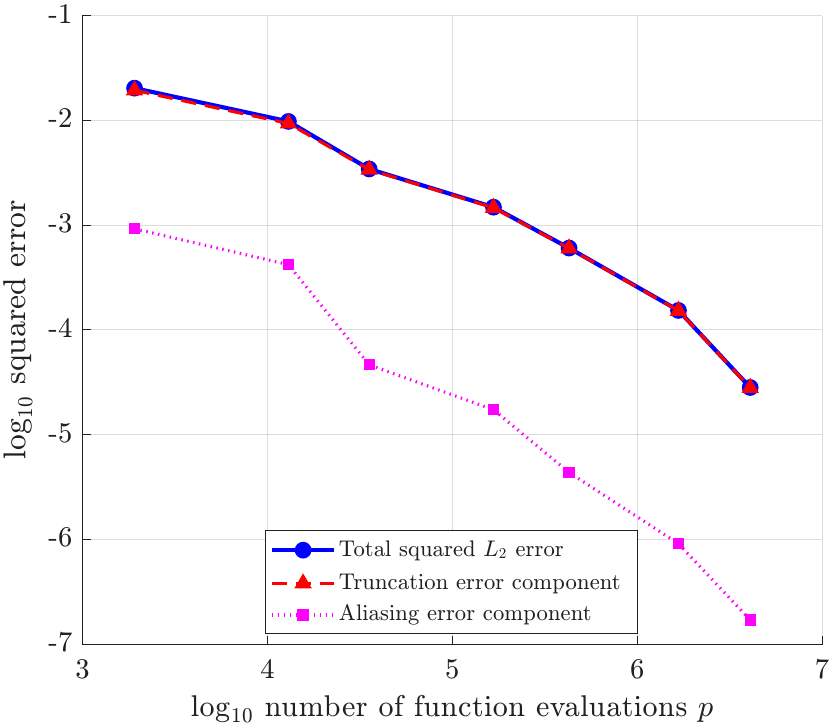}
    \subcaption{$\omega_j=j^{-2}$}
    \label{fig:error_slow}
  \end{minipage}
  \hfill
  \begin{minipage}[b]{0.48\textwidth}
    \centering
    \includegraphics[width=\textwidth]{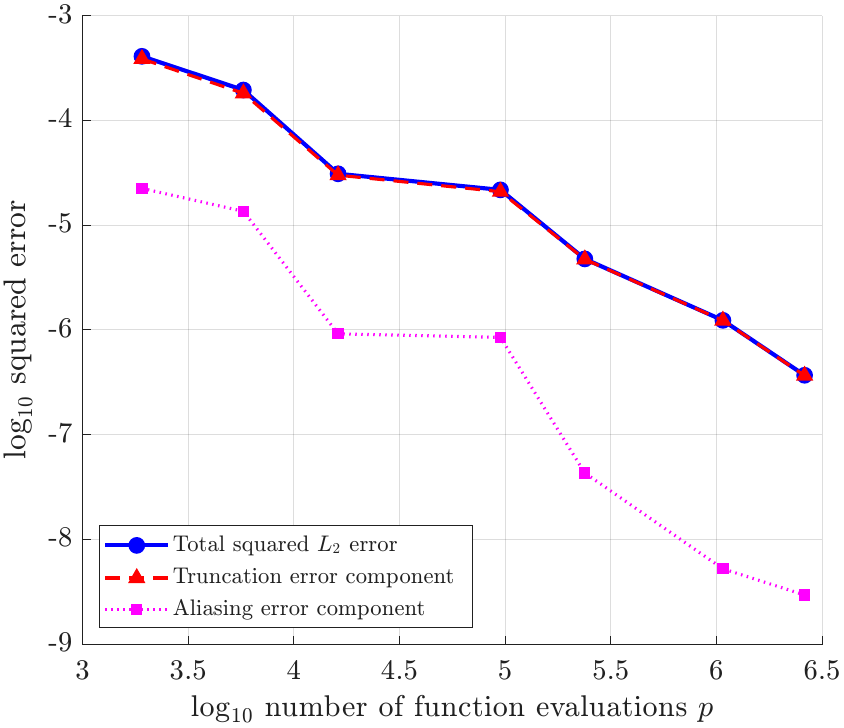}
    \subcaption{$\omega_j=j^{-4}$}
    \label{fig:error_fast}
  \end{minipage}
\caption{The squared $L_2$ error and its aliasing and truncation components for high-dimensional weighted functions with $d=20$.} \label{fig:error_component}
\end{figure}

We further compute the squared $L_2$ error precisely based on the exact Fourier coefficients of $f$. That is, we use the error decomposition
\begin{align*}
    \|f-\mathcal{A}(f)\|_{L_2}^2 & = \sum_{\bsk \in A_{\alpha, \bsgamma,M}} |\widehat{f}(\bsk) - \mathcal{B}_N(f, \bsk)|^2 + \sum_{\bsk \in \mathbb{Z}^d \setminus A_{\alpha, \bsgamma,M}} |\widehat{f}(\bsk)|^2\\
    & = \sum_{\bsk \in A_{\alpha, \bsgamma,M}} |\widehat{f}(\bsk) - \mathcal{B}_N(f, \bsk)|^2 + \left(\|f\|^2_{L_2}-\sum_{\bsk \in A_{\alpha, \bsgamma,M}} |\widehat{f}(\bsk)|^2\right),
\end{align*}
where the first sum on the right-hand side represents the aliasing error, and the remaining terms represent the truncation error. In Figure~\ref{fig:error_component}, we compare the total squared $L_2$ error with its aliasing and truncation components for the case $d=20$. It can be seen that the aliasing error is several orders of magnitude smaller than the truncation error, so that the squared $L_2$ error essentially comes from the truncation error. This small aliasing error demonstrates that our proposed algorithm, utilizing the multiple-shift least-squares approximation, succeeds in recovering the Fourier coefficients within the index set $A_{\alpha, \bsgamma,M}$ with exceptional accuracy. As a result, for the current example, the bottleneck for further error reduction lies not in the recovery process, but rather in expanding the index set to mitigate the dominant truncation error.

\subsection{Computational complexity}
\label{subsec:computational_complexity}

We briefly discuss the computational cost of the proposed algorithms. Let
\[
A := |A_{\alpha,\boldsymbol{\gamma},M}|,\qquad
\Gamma_{\alpha,\boldsymbol{\gamma},N}(\boldsymbol{\ell}^{(j)};\boldsymbol g)
= \{\boldsymbol{\ell}^{(j)}_1, \ldots, \boldsymbol{\ell}^{(j)}_{v_j}\}, \qquad j=1,\ldots,J,
\]
be the partition of \(A_{\alpha,\boldsymbol{\gamma},M}\) into fibers, and let $
R := R_{\alpha,\boldsymbol{\gamma},N}(\boldsymbol g)
= \max_{1\le j\le J} v_j $.
Thus $\sum_{j=1}^J v_j=A$ and $v_j\le R$. In the theoretical construction one can choose a generating vector such that
\[
R=\mathcal{O}((\log N)^{d-1}),
\]
and the repetition parameter is chosen as
\[
S = \left\lceil2 K R \log N\right\rceil
= \mathcal{O}((\log N)^d)
\]
for a fixed constant $K>1$. In the numerical experiments reported above, the observed values of $R$ are small; see the corresponding columns in the tables.

We state the computational cost in the following in terms of the number of function evaluations $p = RSN$.

We separate the cost into three parts. First, the fibers can be constructed by computing the residues
\[
\boldsymbol k\cdot \boldsymbol g \pmod N, \qquad \boldsymbol k\in A_{\alpha,\boldsymbol{\gamma},M},
\]
and grouping equal residues. This requires $\mathcal{O}(dA) \le \mathcal{O}(d p)$ arithmetic operations, up to the cost of sorting or hashing the residues.

Second, one has to compute the quantities
\[
\mathcal{F}_N(f, \boldsymbol{\ell}_m,\boldsymbol y_m^{(s)})
= \frac1N\sum_{n=0}^{N-1}
f\left(\left\{\frac{n\boldsymbol g}{N}+\boldsymbol y_m^{(s)}\right\}\right)
\mathrm{e}^{-2\pi \mathrm{i} n \boldsymbol{\ell}_m \cdot  \boldsymbol{g}/N }.
\]
For fixed $m, s$, one may evaluate \(f\) on the shifted lattice
\[
\left\{\frac{n\boldsymbol g}{N}+\boldsymbol y_m^{(s)}\right\},
\qquad n=0,\ldots,N-1,
\]
and compute all \(N\) residue-class Fourier sums by a one-dimensional FFT of length \(N\). Since the same shifts \({\boldsymbol y_m^{(s)}}\) can be used for all fibers and \(m=1,\ldots,R\), the number of distinct shifted lattices is \(RS\). Therefore the cost of this stage is
$
\mathcal{O}(RSN) = \mathcal{O}(p)
$
function evaluations and
$
\mathcal{O}(RSN\log N) \le \mathcal{O}(p \log p)
$
arithmetic operations for the FFTs with a memory cost of $\mathcal{O}(p)$. 
For the randomized algorithm, if the final error is estimated by a root-mean-square calculation using \(Q\) independent random shifts \(\Delta\), this part of the cost is multiplied by \(Q\).

Third, for each fiber one solves a least-squares problem with matrix
\[
B_j := B(\boldsymbol{\ell}^{(j)},{\boldsymbol y_m^{(s)}})
\in \mathbb{C}^{v_jS\times v_j}.
\]
There is no simple Fourier or Toeplitz structure in these matrices in general, since the frequencies inside one fiber and the shifts are not arranged on a tensor-product grid. Thus these are treated as dense, but relatively small, least-squares problems. If the normal matrices \(B_j^H B_j\) and their factorizations are formed explicitly, the preprocessing cost is bounded by
$
\mathcal{O}\left(
S\sum_{j=1}^J v_j^3
\right)
\le
\mathcal{O}\left(SR^2A\right) \le \mathcal{O}(p R) \le \mathcal{O}(p (\log p)^{d-1})
$, since $A \le N$.
This preprocessing depends only on \(A_{\alpha,\boldsymbol{\gamma},M}\), the generating vector, and the chosen shifts, and can therefore be reused for multiple functions \(f\). Once these factorizations are available, applying the least-squares reconstruction to the right-hand sides costs
$
\mathcal{O}\left(
S\sum_{j=1}^J v_j^2
\right)
\le
\mathcal{O}\left(SRA\right) \le \mathcal{O}(p)
$
arithmetic operations per function, again since $A \le N$. 
If the dense factorizations are included in the one-time setup cost, the worst-case bound is $\mathcal{O}(p (\log p)^{d-1})$.

Hence, after the shifts and least-squares factorizations have been precomputed, the dominant per-function cost is the computation of the shifted lattice Fourier sums, namely
$
\mathcal{O}(RS N\log N) \le \mathcal{O}(p \log p)
$
together with an additional dense least-squares application cost of order
$
\mathcal{O}(SRA) \le \mathcal{O}(p)
$
and a memory requirement of $\mathcal{O}(p)$. 
The least-squares step is therefore a genuine additional cost compared with a standard single-lattice Fourier approximation. 

The preprocessing cost of $\mathcal{O}(p (\log p)^{d-1})$ can be prohibitive in large dimensions, but in the present algorithm the matrices $B_j$ are relatively small because their size is governed by the maximal fiber length $R$, which is polylogarithmic in \(N\) in the theoretical construction and small in the numerical experiments.

Another challenge is that the available bound $p=SRN=\mathcal{O}(N(\log N)^{2d-1})$ has exponential dependence on the dimension $d$.

\section{Conclusion and outlook}

We have introduced a new lattice-based method for the approximation of periodic functions in weighted Korobov spaces. The method uses shifted copies of a single rank-1 lattice and resolves the associated aliasing by partitioning the frequency index set into fibers and recovering the Fourier coefficients within each fiber from a least-squares problem. The resulting deterministic and randomized algorithms attain, up to logarithmic factors and an arbitrarily small loss in the exponent, the optimal polynomial convergence orders for $L_\infty$- and randomized $L_2$-approximation, respectively. Thus, the present work shows that a single underlying lattice, when combined with suitably chosen shifts, can overcome limitations that arise for approximation based on an unshifted rank-1 lattice.

The theoretical results should, however, be viewed together with several practical limitations. Most importantly, the currently available bound for the number of function evaluations is $p = N R_{\alpha,\boldsymbol{\gamma},N}(\boldsymbol{g})S = \mathcal{O}\bigl(N(\log N)^{2d-1}\bigr)$. For a fixed dimension, this represents a polylogarithmic sampling overhead relative to the lattice size $N$. Nevertheless, the exponent of the logarithmic factor grows linearly with $d$, so the worst-case sampling requirement deteriorates rapidly as the dimension increases. For comparison, the oversampling factors established in \cite{Kammerer2018JFAA, BGKS25} are of logarithmic order. Hence our method in its present form is most directly applicable to problems of moderate dimension. 

A second practical issue is the cost of recovering the coefficients within the individual fibers. The size of the least-squares systems is determined by the fiber lengths and by the number of shifts. Although the numerical experiments indicate that the observed fiber lengths are often considerably smaller than the general theoretical upper bounds, large fibers can increase both the computational work and the memory requirements. Moreover, the use of normal equations, if implemented directly, may be less numerically stable than solving the corresponding least-squares problems using QR or other stable factorizations. The evaluations on different shifts and the computations associated with different fibers are largely independent, however, and therefore offer natural opportunities for parallel implementation.

These observations suggest several possible algorithmic improvements. One direction is to find explicit constructions of shifts $\bm{y}_m^{(s)}$, thereby obtaining a fully deterministic approximation algorithm.

To reduce the preprocessing cost, one could use an iterative least-squares solver to reduce the cost from $\mathcal{O}(p R)$ to $\mathcal{O}(p k)$ using $k$ iterations. It would also be useful to develop sharper probabilistic estimates for the typical fiber lengths, since the numerical results suggest that worst-case bounds may substantially overestimate the cost encountered in practice.

The extension to genuinely high-dimensional settings remains an important open problem. The use of weighted Korobov spaces already provides a mechanism for describing decreasing importance of successive variables. Future work could investigate whether assumptions on the decay of the weights lead to dimension-independent or weakly dimension-dependent bounds for the fiber lengths, the required number of shifts, and the overall sampling complexity. Dimension-adaptive frequency sets, product-and-order-dependent weights, or constructions based on effective dimension may make the approach applicable when the nominal dimension is large. Establishing tractability results for the proposed algorithms would clarify precisely under which weight conditions the dimension dependence can be controlled.

Further extensions include the treatment of other periodic smoothness spaces and the use of periodisation or change-of-variable techniques for nonperiodic approximation problems. 

Overall, the results provide a theoretical foundation for Fourier approximation using multiple shifts of a single rank-1 lattice, while the reduction of the sampling overhead, dimension dependence, and coefficient-recovery costs remains central to making the method scalable for larger and higher-dimensional applications.

\bibliography{ref}                        


\end{document}